\documentclass[a4paper,english]{amsart}

\usepackage{amsfonts,amsmath,amssymb,color}
\usepackage{enumerate}

\title{Sobolev Inequalities for Differential Forms and  $L_{q,p}$-cohomology.}
\date{June 3, 2005 modified July 19, 2006} 

\author{Vladimir Gol'd'shtein}
\address{Department of Mathematics,
Ben Gurion University of the Negev,
P.O.Box 653, Beer Sheva, Israel} 
\email{vladimir@bgumail.bgu.ac.il}

\author{Marc Troyanov}
\address{Marc Troyanov     
Institut de G{\'e}om{\'e}trie, alg{\`e}bre et topologie (IGAT)
B{\^a}timent BCH \\
 \'Ecole Polytechnique F{\'e}derale de
Lausanne, 1015 Lausanne - Switzerland}
\email{marc.troyanov@epfl.ch}


 \newtheorem{thm}{Theorem}[section]
 \newtheorem{defn}[thm]{Definition}
 \newtheorem{rem}[thm]{Remark}
 \newtheorem{prop}[thm]{Proposition} 
 \newtheorem{lem}[thm]{Lemma} 
 \newtheorem{cor}[thm]{Corollary} 

 \numberwithin{equation}{section}
 \renewcommand{\r}{\mathbb{R}}
 \newcommand{\ds}{\displaystyle}
 \DeclareMathOperator{\Imm}{Im}
 \DeclareMathOperator{\vol}{vol}

\makeatother

\begin{document}

\begin{abstract}
We study the relation  between  Sobolev inequalities for differential forms on  a Riemannian manifold  $(M,g)$
and the $L_{q,p}$-cohomology of  that manifold. \\
The $L_{q,p}$-cohomology of  $(M,g)$ is defined to be the quotient of the space
of closed differential forms in $L^p(M)$ modulo the exact forms which are exterior differentials of forms in 
$L^q(M)$. 
\end{abstract}

\maketitle


\section{Introduction}

Let us start by stating a  Sobolev type Inequality for differential forms on a compact manifold:

\begin{thm}\label{th.sobin1}
Let $(M,g)$ be a smooth $n$-dimensional compact Riemannian manifold,
$1\leq k \leq n$ and $p,q \in (1,\infty)$. Then there exists a constant $C$
such that for any  differential form $\theta$ of degree
$k-1$ on $M$ with coefficients in $L^{q}$, we have
\begin{equation}\label{ineq.sob1}
  \inf_{\zeta \in Z^{k-1}} \|\theta-\zeta\|_{L^q(M)} \leq C  \|d\theta\|_{L^p(M)},
\end{equation}
if and only if
\begin{equation} \label{cond.sobolev}
  \frac{1}{p}-\frac{1}{q} \leq \frac{1}{n}.
\end{equation}
Here $Z^{k-1}$ denotes the set of smooth closed $(k-1)$-forms on $M$.
\end{thm}

The differential $d\theta$ in the inequality above is to be understood in the sense of currents.

\medskip

Note that condition (\ref{cond.sobolev}) is equivalent to
\begin{equation} \label{cond.sobolev+}
  p\geq n \qquad \text{or} \qquad p < n  \ \text{ and } \ q \leq
  p^*=\frac{np}{n-p}.
\end{equation}

In the case of zero forms (i.e. $k=1$), this theorem can be
deduced from the corresponding result for functions with compact
support in $\r^n$ by a simple argument using a partition of unity.
The case of differential forms of higher degree can be proved
using more involved reasoning based on standard results from the
Hodge--De Rham theory and $L^p$-elliptic estimates obtained in
the 1950' by various authors. We give a sketch of such a proof in
the appendix of this paper.

\medskip

In the case of a non compact manifold, the inequality
(\ref{ineq.sob1}) is still meaningful if the differential form
$\theta$ belongs to $L^q$. Although the condition
(\ref{cond.sobolev}) is still necessary in the non compact case,
it is no longer sufficient and additional conditions must be
imposed on the geometry of the manifold $(M,g)$ for a Sobolev
inequality to hold.

\medskip

The main goal of this paper is to investigate these
conditions. Our Theorem \ref{th.sobin2a} below gives a necessary and
sufficient condition based on an invariant called the
\emph{$L_{q,p}$--cohomology of $(M,g)$} and which is defined as
\begin{equation*}
  H^k_{q,p}(M) = Z^k_{p}(M)/d\Omega^{k-1}_{q,p}(M).
\end{equation*}
where $Z^k_{p}(M)$ is the Banach space of closed $k$-forms
$\theta$ in $L^p(M)$ and $\Omega^{k-1}_{q,p}(M)$ is the space of
all $(k-1)$-forms $\phi$ in $L^q(M)$ such that $d\phi\in L^p$.

\medskip

We will also prove a regularization theorem saying that any
$L_{q,p}$-cohomology class can be represented by a smooth form,
provided that (\ref{cond.sobolev}) holds (see Theorem
\ref{th.smoothcohomology}).  This implies in particular that the
$L_{q,p}$-cohomology of a compact manifold $M$ coincides with the
usual De Rham cohomology $M$ and it gives us a new proof of
Theorem \ref{th.sobin1} above. This new proof is perhaps simpler
than the classical one sketched in the appendix (at least it does
not rely on the rather deep elliptic estimate).

\medskip

The techniques of this paper also provide a proof of the
following result which is a complement to Theorem \ref{th.sobin1}:

\begin{thm}\label{th.sobin2} \  Let $(M,g)$ be a smooth compact Riemannian manifold of dimension
$n$ and $p,q \in (1,\infty)$.   There exists a constant $C$ such that for all closed
differential forms $\omega$ of degree $k$ with coefficients in $L^{p}(M)$, there exists a
differential form $\theta $ of degree $k-1$ such that $d\theta =
\omega$ and
\begin{equation}\label{ineq.sob2}
 \left\| \theta \right\| _{L^{q}}\leq C\left\| \omega
 \right\|_{L^{p}},
\end{equation} if and only if $p,q$ satisfy the condition
(\ref{cond.sobolev}) and
$H_{\text{\tiny{\emph{DeRham}}}}^{k}(M)=0$.
\end{thm}

\medskip

Both  Theorems \ref{th.sobin1} and  \ref{th.sobin2} are proved at the end of section \ref{sec.regularization}.
In the non compact case, we  prove in Theorem \ref{th.sobin1a} below that the inequality (\ref{ineq.sob2})  holds if and only if
$H_{q,p}^{k}(M,g)=0$.

\bigskip

The Sobolev inequality is important because it is a key
ingredient in solving  partial differential equations. To
illustrate this point, we show in  section \ref{sec.PDE} how
Theorem \ref{th.sobin2a} can be used to solve the non linear
equation
\begin{equation}\label{smlapl}
 \delta(\left\Vert d\theta\right\Vert ^{p-2}d\theta)=\alpha
\end{equation}
for differential forms. Here $\delta$ is the formal adjoint to the exterior differential $d$.

\medskip

Although it is certainly a nice observation that such Sobolev
type inequalities for differential forms have interpretations in
$L_{q,p}$-cohomology, this will not lead us very far unless we are
able to compute some of this cohomology. Unfortunately, this is
not an easy task and only few examples of $L_{q,p}$-cohomology
groups are presently known. It is thus also one of our goals in
this paper to begin developing some of the basic facts from the
theory. In particular, we present here some results in the direction of
duality (see section \ref{sec.almostduality}), a proof of the
Poincar\'e Lemma for $L_{q,p}$-cohomology and a non vanishing
result for the $L_{q,p}$-cohomology of the hyperbolic plane
$\mathbb{H}^2$. This non vanishing result says in particular that
the Sobolev inequality (\ref{ineq.sob2}) for one-forms  never
holds on $\mathbb{H}^2$ for any $p,q\in (1,\infty)$.

\bigskip

{\footnotesize 
\begin{enumerate}
\item[\textbf{Content}]   
\item[\textmd{1}.]  {Introduction} 
\item[\textmd{2}.]  {Definitions} 
\item[\textmd{3}.] {Some elementary properties of $L_{q,p}$-cohomology} 
\item[\textmd{4}.]  {Banach complexes}
\item[\textmd{5}.]  {$L_{q,p}$-cohomology and Banach complexes} 
\item[\textmd{6}.]  {$L_{q,p}$-cohomology and Sobolev inequality} 
\item[\textmd{7}.]  {Manifolds with finite volume and monotonicity} 
\item[\textmd{8}.]  {Almost duality}
\item[\textmd{9}.]  {The $L_{q,p}$-cohomology of the line}
\item[\textmd{10}.]  {The cohomology of the hyperbolic plane} 
\item[\textmd{11}.]  {The cohomology of the ball} 
\item[\textmd{12}.]  {Regularization of forms and cohomology classes}
\item[\textmd{13}.]  {Relation with a non linear PDE}
\item[\textmd{14}.] {Torsion in $L_{2}$-cohomology and the Hodge-Kodaira decomposition}
\item[\textmd{Appendix}.]  { A ``classic'' proof of Theorem 1.1 in the compact case.}
\end{enumerate}
}

\medskip

Let us shortly describe what is contained in the paper.
In sections 2 and 3, we give the necessary definitions and we prove some
elementary properties of $L_{q,p}$-cohomology. Then we present
some basic facts of the theory of Banach complexes and we derive
the cohomological interpretation of Sobolev inequalities for
differential forms (section 4,5 and 6). In section 7, we prove
some monotonicity properties for the $L_{q,p}$-cohomology of
finite dimensional manifolds and in section 8 we introduce a
notion of  ``almost duality'' techniques (a standard Poincar\'e
duality holds only when  $p=q$). We apply these techniques to
compute the $L_{q,p}$-cohomology of the line (section 9) and the hyperbolic plane (section
10) and to prove a version of the Poincar\'e Lemma (section 11).
In section 12, we show that the $L_{q,p}$-cohomology of a manifold
can be represented by smooth forms under the condition (\ref{cond.sobolev}).
Finally, we show in section 13 how the $L_{q,p}$-cohomology can
be relevant in the study of some non linear PDE, and in section 14 we
give a relation between the $L_2$-cohomology and the Laplacian on
complete manifolds.
The paper ends with an appendix describing  an
alternative proof of Theorems \ref{th.sobin1} based on $L^p$
elliptic estimates.

\medskip

\textbf{Remark.}  The reader might prefer to call the inequality (\ref{ineq.sob1}) a \emph{Poincar\'{e} inequality} and use the term \emph{Sobolev inequality} only for the inequality  (\ref{ineq.sob2}). In fact there are various uses of the terms Poincar\'{e} and Sobolev inequalities. According to \cite{GT},
the Poincar\'{e} inequality is simply a special case of the Sobolev one (it is
in fact the case $p=q$). In this paper, we avoid the  name Poincar\'{e} inequality.

\bigskip

\emph{Acknowledgment.}  \ Part of this research has been done in
the autumn of 2001, when both authors stayed at IHES in
Bures-Sur-Yvette. We are happy to thank the Institute for its warm
hospitality. We also thank Pierre Pansu for his interest in our
work and for the kindness and patience with which he explained us
his viewpoint on the subject.

\section{Definitions}

Let us recall the notion of weak exterior differential of a
differential form on a Riemannian manifold $(M,g)$.

We denote by $C^{\infty}_{c}(M,\Lambda^{k})$ the vector space
of smooth differential forms of degree $k$ with compact support on
$M$ and by $L^1_{loc}(M,\Lambda^k)$ the space of differential
$k$-forms whose coefficients (in any local coordinate system) are
locally integrable.

\medskip

\begin{defn}
One says that a form  $\theta \in L^1_{loc}(M,\Lambda^k)$ is the
\emph{weak exterior differential} of a form  $\phi \in
L^1_{loc}(M,\Lambda^{k-1})$ and one writes $d\phi = \theta$ if for
each $\omega \in C^{\infty}_{c}(M,\Lambda^{n-k})$, one has
\[
 \int_M \theta \wedge \omega = (-1)^{k}\int_M \phi \wedge d\omega
 \,  .
\]
\end{defn}

Clearly $d\phi$ is uniquely determined up to sets of Lebesgue measure zero, because
$d\phi$ is the exterior differential (in the sense of currents) of the current $\phi$.
It is also clear that $d\circ d =0$, and this fact  allows us to define
various cohomology groups.

\medskip

Let  $L^p(M,\Lambda^k)$ be the space of differential forms in
$L^1_{loc}(M,\Lambda^k)$ such that
\[
 \| \theta \|_p := \left( \int_M |\theta |^pdx\right)^{\frac{1}{p}}
  < \infty \, .
\]
We then set  $Z^k_p(M) :=L^p(M,\Lambda^k) \cap \ker d$ (= the set
of weakly closed forms in $L^p(M,\Lambda^k)$) and
\[
  B^k_{q,p}(M) := d\left(L^{q}(M,\Lambda^{k-1}) \right) \cap
  L^p(M,\Lambda^k).
\]

\begin{lem}
  $Z^k_p(M) \subset L^p(M,\Lambda^k)$ is a closed linear subspace. 
  In particular it is a Banach space.
\end{lem}

\textbf{Proof}  We need to show that an arbitrary element $\phi \in \overline{Z}^k_p(M)$ in the closure of 
$Z^k_p(M)$ is a weakly closed form. Choose a sequence $\phi_{i} \in {Z}^k_p(M)$ such that $\phi_{i} \to \phi$ in $L^p$-norm.
Since $\phi_{i}$ are weakly closed forms, we have
\[
  \int_M \phi_{i} \wedge d\omega = 0,
\]
for any smooth differential forms $\omega$ of degree $n-k-1$ with
compact support on $M$. Using H\"older's inequality, we obtain
\[
  \int_M \phi \wedge d\omega  =   \int_M (\phi -\phi_{i}) \wedge d\omega 
  \leq \|\phi -\phi_{i} \|_{L^p(M)} \|d\omega \|_{L^{p'}(M)}
  \to 0.
\]
Here $1/p+1/p'=1$.

Thus $ \int_M \phi \wedge d\omega  = 0$ for any 
$\omega = C^{\infty}_{c}(M,\Lambda^{n-k-1})$   and hence $\phi \in {Z}^k_p(M)$.

\qed

\medskip

Observe that   $B^k_{q,p}(M) \subset Z^k_p(M)$ (because $d\circ d = 0$), we thus have
 $$B^k_{q,p}(M) \subset \overline{B}^k_{q,p}(M) \subset Z^k_p(M)
 =  \overline{Z}^k_p(M) \subset L^p(M,\Lambda^k).$$
 
\medskip

\begin{defn}
The $L_{q,p}$\emph{-cohomology} of $(M,g)$ (where $1\leq p,q\leq
\infty $) is defined to be the quotient
\begin{equation*}
H_{q,p}^{k}(M):=Z_{p}^{k}(M)/B_{q,p}^{k}(M)\,,
\end{equation*}
and the \emph{reduced $L_{q,p}$-cohomology} of $(M,g)$ is
\begin{equation*}
\overline{H}_{q,p}^{k}(M):=Z_{p}^{k}(M)/\overline{B}_{q,p}^{k}(M)\,,
\end{equation*}
(where $\overline{B}_{q,p}^{k}(M)$ is the closure of
$B_{q,p}^{k}(M)$). We also define the \emph{torsion} as
$$T_{q,p}^{k}(M):=\overline{B}_{q,p}^{k}(M)\,/ B_{q,p}^{k}(M).
$$
\end{defn}
We thus have the exact sequence
\begin{equation*}
0\rightarrow T_{q,p}^{k}(M)\rightarrow H_{q,p}^{k}(M)\rightarrow
\overline{H}_{q,p}^{k}(M)\rightarrow 0.
\end{equation*}

\bigskip

The reduced cohomology is naturally a Banach space. The unreduced
cohomology is a Banach space if and only if the torsion vanishes.

\medskip

By Lemma \ref{lem.notorsion} below, we see that the 
torsion $T_{q,p}^{k}(M)$ can be either $\{0\}$ or
infinite dimensional. Indeed, if $\dim T_{q,p}^{k}(M)<\infty$ then
$B_{q,p}^{k}(M)$ is closed, hence $T_{q,p}^{k}(M)=\{0\}$. In
particular, if $\dim T_{q,p}^{k}(M)\neq 0$ then $\dim
H_{q,p}^{k}(M)=\infty$.

\medskip

When $p = q$, we simply speak of $L_{p}$-cohomology and write
$H^k_{p}(M)$ and $\overline{H}^k_{p}(M)$.

\bigskip

\textbf{Example} The  $L_{q,p}$-cohomology of the bounded interval
$M=(0,1)$ is easily computed: we clearly have $H_{q,p}^{0}((0,1))=\mathbb{R}$ and 
  $H_{q,p}^{1}((0,1))=0$ for any $1\leq q,p\leq\infty$.

Indeed if  $\omega=a(x)dx$ belongs to  $L^{p}((0,1))\subset L^{1}((0,1))$,
then  $f(x):=\int_{-\infty}^{x}a(s)ds$
belongs to $L^{q}((0,1))$ for any  $1\leq q\leq\infty$.

\medskip 

The   $L_{q,p}$-cohomology of the unbounded intervals and other examples will be computed below.

\section{Some elementary properties of $L_{q,p}$-cohomology}

\subsection{Zero dimensional cohomology.}

We have
$H_{q,p}^{0}(M)=\overline{H}_{q,p}^{0}(M)=Z_{p}^{0}(M)=H_{p}^{0}(M)$
and these spaces have the following interpretation: \ $\dim
H_{\infty}^{0}(M)$ is the number of connected components of $M$
and $\dim H_{p}^{0}(M)$ is the number of connected components
with finite volume of $M$ if $1\leq p<\infty$.

\subsection{Conformal invariance.}

Let $(M,g)$ be a Riemannian manifold of dimension $n$. Recall that a new metric
$g_{1}$ is a conformal deformation of $g$ if $g_{1}:=\rho^{2}g$
where $\rho:M\rightarrow\mathbb{R}_{+}$ is a smooth function.

\medskip

The pointwise norms of a $k$-form $\omega$ with respect to the metrics
$g_{1}$ and $g$ are related by the identity $|\omega|_{g_{1}}=\rho^{-k}|\omega|_{g}$.
The volume elements are related by \ $d\vol_{g_{1}}=\rho^{n}d\vol_{g}$.
In particular
\[
|\omega|_{g_{1}}^{p}d\vol_{g_{1}}=\rho^{n-pk}|\omega|_{g}^{p}d\vol_{g}
\]
for any $k$-form; likewise,
$|\theta|_{g_{1}}^{q}d\vol_{g_{1}}=\rho^{n-q(k-1)}|\theta|_{g}^{q}d\vol_{g}$
for any $k-1$-form $\theta$. It follows that
$H_{q,p}^{k}(M,g_{1})=H_{q,p}^{k}(M,g)$ if $n-pk=n-q(k-1)=0$.

\medskip

We thus have the

\begin{thm}\label{th:conf}
If  $q=\frac{n}{k-1}$ and $p=\frac{n}{k}$, then 
$H_{q,p}^{k}(M,g)$ and $\overline{H}_{q,p}^{k}(M,g)$ are conformal invariants.
\end{thm}

\qed

\section{Banach complexes}

The abstract theory of Banach complexes is based on a combination
of techniques from homological algebra and functional analysis;
this theory is the natural framework of $L_{q,p}$-cohomology and
we shall take this point of view to show the connections between
Sobolev inequalities and $L_{q,p}$-cohomology.

\medskip

There is not much literature on Banach complexes, we therefore
give below all necessary definitions. The reader may look in
\cite{Grom} for more information.

\subsection{Cohomology of Banach complexes and abstract Sobolev inequalities.}

\begin{defn}\label{def:ban.comp}
A \emph{Banach complex} is a sequence $F^{*}=\{ F^{k},d_{k}\}_{k\in\mathbb{N}}$ \
where $F^{k}$ is a Banach space, $d_{k}:F^{k}\rightarrow F^{k+1}$ is a
bounded operator and $d_{k+1}\circ d_{k}=0$.
\end{defn}

\textbf{Remarks} 1.) It would be more correct to call such an object
a Banach cocomplex (and to use the name complex for the case where
$d_{k}:F^{k}\rightarrow F^{k-1}$), but for simplicity, we shall speak of
complexes.

2) To simplify notations, we usually note $d$ for any of the
operators $d_{k}$.

\bigskip

\begin{defn}
Given a Banach complex $\{ F^{k},d\}$ \ we introduce the following vector spaces:
\begin{itemize}
\item $Z^{k}:=\ker(d:F^{k}\rightarrow F^{k+1})$, it is a closed subspace
of $F^{k}$;
\item $B^{k}:=$Im$(d:F^{k-1}\rightarrow F^{k})\subset Z^{k}$;
\item $H^{k}(F^{*}):=Z^{k}/B^{k}$ is the \emph{cohomology} of the complex
$F^{*}=\{ F^{k},d\}$;
\item $\overline{H}^{k}(F^{*}):=Z^{k}/\overline{B}^{k}$ is \emph{the reduced
cohomology} of the complex $F^{*}$;
\item $T^{k}(F^{*}):=\overline{B}^{k}/B^{k}=H^{k}/\overline{H}^{k}$ is
the  \emph{torsion} of the complex $F^{*}$.
\end{itemize}
\end{defn}

\medskip

Let us make a few elementary observations :

\begin{enumerate}[a.)]
\item $\overline{H}^{k},Z^{k}$ and $\overline{B}^{k}$ are Banach spaces;
\item The natural (quotient) topology on $T^{k}:=\overline{B}^{k}/B^{k}$
is coarse (any closed set is either empty or $T^{k}$);
\item We have the exact sequence
\begin{equation*}\label{exact}
  0\rightarrow T^{k}\rightarrow
H^{k}\rightarrow\overline{H}^{k}\rightarrow 0.
\end{equation*}
\end{enumerate}

\bigskip

There is a natural notion of subcomplex:
\begin{defn}
 A \emph{subcomplex} $G^*$ of a Banach complex $\{F^*,d\}$ is a sequence  
 of linear subspaces $G^k\subset F^k$ (not necessarily closed) such
 that $d(G^k)\subset G^{k+1}$.  If all 
$G^k$ are closed subspaces, we say that $G^*$ is a 
\emph{Banach-subcomplex} of $F^*$.\\
\end{defn}

\medskip

 The cohomology of the subcomplex $G^*$ is defined as
 $$H^k(G^*)=(G^k\cap\ker d)/d(G^{k-1}).$$

\medskip

Observe that in general $H^k(G^*)$ is not a Banach space, but
there is no way to define a reduced cohomology of $G^*$, unless
$G^*\subset F^*$ is a Banach-subcomplex.

\medskip

\begin{lem}\label{lem.notorsion}
For any Banach complex $\{ F^{k},d\}$, the following conditions
are equivalent
\begin{enumerate}[(i.)]
\item $T^{k}=0$;
\item $\dim T_{k} < \infty$;
\item $B^{k}  \subset F^{k}$ is closed.
\end{enumerate}
\end{lem}

\textbf{Proof} (i)$\Rightarrow$(ii) is obvious and
(ii)$\Rightarrow$(iii) follows  e.g  from \cite[Th. 3.2 page 27]{Edmunds87}.
The implication  \ (iii)$\Rightarrow$(i) follows directly from the definition of the torsion.

\qed

\bigskip

\begin{prop}\label{abs.sob1}
The following are equivalent:
\begin{enumerate}[(i)]
  \item $H^{k}=0$;
  \item The operator $d_{k-1}:F^{k-1}/Z^{k-1}\rightarrow Z^{k}$ admits a
bounded inverse $d_{k-1}^{-1}$;
  \item There exists a constant $C_{k}$ such that for any \ $\theta\in
Z^{k}$ there is an element $\eta\in F^{k-1}$ with $d\eta=\theta$
and
\[
 \Vert\eta\Vert_{F^{k-1}}\leq C_{k}\Vert\theta\Vert_{F^{k}}.
\]
\end{enumerate}
\end{prop}

\textbf{Proof}
(i) $\Rightarrow$ (ii). \ Suppose $H^{k}=0$. Then
$d_{k-1}:F^{k-1}/Z^{k-1}\rightarrow Z^{k}$ is a bijective bounded
linear operator and by the open mapping theorem, the inverse map
\[ d_{k-1}^{-1}:Z^{k}\rightarrow F^{k-1}/Z^{k-1}\]
is also a bounded operator.

\medskip

(ii) $\Rightarrow$ (iii). \ Let $\gamma$ be the norm of
$d_{k-1}^{-1}:Z^{k}\rightarrow F^{k-1}/Z^{k-1}$, then for any
$\theta\in Z^{k}$ we can find $\xi\in F^{k-1}$ such that
$d_{k-1}\xi=\theta$. Furthermore
\[
\left\Vert [\xi]\right\Vert
_{F^{k-1}/Z^{k-1}}=\inf_{\zeta\in Z^{k-1}}\left\Vert
\xi-\zeta\right\Vert _{F^{k-1}}\leq\gamma\left\Vert
\theta\right\Vert _{F^{k}}.
\]
In particular, there exists $\zeta\in Z^{k-1}$ such that
$\left\Vert \xi-\zeta\right\Vert _{F^{k-1}}\leq 2\gamma\left\Vert
\theta\right\Vert _{F^{k}}$. Let us set $\eta:=(\xi-\zeta)$, then
$d_{k-1}\eta=\theta$ and $\left\Vert \eta\right\Vert
_{F^{k-1}}\leq C_{k}\left\Vert \theta\right\Vert _{F^{k}}$with
$C_{k}=2\gamma=2\left\Vert d_{k-1}^{-1}\right\Vert
_{Z^{k}\rightarrow F^{k-1}/Z^{k-1}}$.

\medskip

The implication (iii) $\Rightarrow$ (i) is clear.

\qed

\medskip

\begin{prop}\label{abs.sob2} The following conditions are equivalent:
\begin{enumerate}[(i)]
  \item $T^{k}=0$;
  \item The operator $d_{k-1}:F^{k-1}/Z^{k-1}\rightarrow B^{k}$ admits a
bounded inverse $d_{k-1}^{-1}$.
\end{enumerate}
And any one of these conditions imply
\begin{enumerate}[(i)] \setcounter{enumi}{2}
  \item There exists a constant $C_{k}^{'}$ such that for any \ $\xi\in F^{k-1}$
 there is an element $\zeta\in Z^{k-1}$ such that
\begin{equation}
\Vert\xi-\zeta\Vert_{F^{k-1}}\leq C_{k}^{'}\Vert
d\xi\Vert_{F^{k}}.\label{eq2}
\end{equation}
\end{enumerate}
\end{prop}

\medskip

\textbf{Proof} The conditions (i) and (ii) are equivalent, because
the existence of a bounded inverse operator is equivalent to the
 closedness of $B^{k-1}$ by the open mapping theorem.

\medskip

Let us assume that $T^{k}=0$ and prove (iii). By hypothesis,
$B^{k}$ is a Banach space and $d_{k-1}:F^{k-1}/Z^{k-1}\rightarrow
B^{k}$ is a bijective bounded linear operator. Thus, by the open
mapping theorem, the inverse $d_{k-1}^{-1}:B^{k}\rightarrow
F^{k-1}/Z^{k-1}$ is also a bounded operator.

Let $\gamma$ be the norm of $d_{k-1}^{-1}:B^{k}\rightarrow
F^{k-1}/Z^{k-1}$, then for any $\xi\in F^{k-1}$ we have
\[
\left\Vert [\xi]\right\Vert _{F^{k-1}/Z^{k-1}}=\inf_{\zeta\in
Z^{k-1}}\left\Vert \xi-\zeta\right\Vert
_{F^{k-1}}\leq\gamma\left\Vert d_{k-1}\xi\right\Vert _{F^{k}}
\]
in particular, there exists $\zeta\in Z^{k-1}$ such that
$\left\Vert \xi-\zeta\right\Vert _{F^{k-1}}\leq2\gamma\left\Vert
d_{k-1}\xi\right\Vert _{F^{k}}$.

\qed

\begin{prop} \label{abs.sob3}
 If  $F^{k-1}$ is a reflexive Banach space, then
the three conditions of the previous proposition are equivalent.
\end{prop}

\medskip

\textbf{Proof}
We only need to show that (iii)$\Rightarrow$ (i) i.e.
$B^{k}=\overline{B}^{k}\subset F^{k}$ provided  (\ref{eq2}) holds
and $F^{k-1}$ is a reflexive. Let $\theta\in\overline{B}^{k}$,
then there exists a sequence $\xi_{i}\in F^{k-1}$ such that
$d_{k-1}\xi_{i}\rightarrow\theta$ in $F^{k}$. By hypothesis there 
exists a sequence $\zeta_{i}\in Z^{k-1}$ such that $\left\Vert
\xi_{i}-\zeta_{i}\right\Vert _{F^{k-1}}\leq C_{k}'\left\Vert
d\xi_{i}\right\Vert _{F^{k}}$. In particular, the sequence
$\{\eta_{i}:=(\xi_{i}-\zeta_{i})\}$ is bounded, we may thus find
a subsequence (still denoted $\{\eta_{i}\}$) which converges
weakly to an element $\eta\in F^{k-1}$.

Using the Mazur Lemma (see e.g. chap. V \S 1,  Theorem 2, page 120 in \cite{yos}), we may construct a sequence
$\{\widetilde{\eta}_{i}=\sum_{j=i}^{N(i)}a_{i}\eta_{j}\}$
of convex combinations of $\eta_{i}$ such that $\widetilde{\eta}_{i}$
converges strongly to $\eta$. We then have
\[
d_{k-1}\eta=\lim_{i\rightarrow\infty}d_{k-1}\widetilde{\eta}_{i}=
\lim_{i\rightarrow\infty}\sum_{j=i}^{N(i)}a_{i}d_{k-1}\eta_{i}=
\lim_{i\rightarrow\infty}\sum_{j=i}^{N(i)}a_{i}d_{k-1}\xi_{j}=\theta
\]
hence $\theta\in\text{Im}(d)=B^{k}$. We proved that $B^{k}$ is
closed, i.e. $T^{k}=0$.
 
\qed

\subsection{Morphisms and  homotopies of Banach complexes.}

This part will be useful to regularize $L_{q,p}$-cohomology, see section \ref{sec.regularization}.

\medskip

\textbf{Definitions 1)} A \emph{morphism}  $R^{\ast}$ between two
Banach complexes $F^{\ast}=\{ F^{k},d\}$ and $E^{\ast}=\{ E^{k},d\}$
is a family of bounded operators $R^{k}:F^{k}\rightarrow E^{k}$ such
that
 \[
 d_{k}\circ R^{k}=R^{k+1}\circ d_{k}.
\]

\textbf{2)} A \emph{homotopy} between two morphisms $R^{\ast}$
and $S^{\ast}:F^{\ast}\rightarrow E^{\ast}$ is a family of bounded
operators $A^{k}:F^{k}\rightarrow E^{k-1}$ such that \[
S^{k}-R^{k}=d_{k-1}\circ A^{k}+A^{k+1}\circ d_{k}.\]

\textbf{3)} A \emph{weak homotopy} between two morphisms
$R^{\ast}$ and $S^{\ast}:F^{\ast}\rightarrow E^{\ast}$ is a
sequence of families of bounded operators
$A_{j}^{k}:F^{k}\rightarrow E^{k-1}$ such that for any element
$x\in F^{k}$ we have
\[
\lim_{j\rightarrow\infty}\left\Vert (d_{k-1}\circ A_{j}^{k}+A_{j}^{k+1}
\circ d_{k})x-(S^{k}-R^{k})x\right\Vert =0.
\]

\medskip

Observe that, if $R^* = \{R^{k}:F^{k}\rightarrow E^{k}\}$ is a morphism, then 
its image is a subcomplex of $E^*$ and it is a  Banach-subcomplex if and
only if all $R^k$ are closed operators. The kernel of $R^*$ is  always a Banach-subcomplex of $F^*$.

\medskip

\begin{prop}\label{prop:chain.com} %
Let $R^{\ast}:F^{\ast}\rightarrow F^{\ast}$ be an endomorphism of
a Banach complex $\{F^*,d\}$ such $R^{\ast}(F^{\ast})\subset G^*$
where $G^*$ is a subcomplex.

\smallskip

If there exists a homotopy $\{A^k : F^k\to F^{k-1}\}$ between
$R^*$ and the identity operator $I : F^* \to F^*$, then
$$H^k(F^*) = H^k(G^*).$$
\end{prop}

\medskip

\textbf{Proof} Given $\xi\in Z^k(F^*)$, we observe that $R^k\xi
\in Z^k(G^*)$ because $dR\xi=Rd\xi=0$. If $\xi=d\eta\in B^k(F^*)$,
then $R^k\xi=R^kd\eta=dR^k\eta\in B^k(G^*)$.

This proves that $[R\xi]$ is a well defined cohomology class in
$H^k(G^*)$ for any cohomology class $[\xi]\in H^k(F^*)$.

But since
$$\xi - R\xi = dA\xi + Ad\xi = dA\xi$$
for any $\xi\in Z^k(F^*)$, we see that in fact $[R\xi]=[\xi]\in
H^k(F^*)$ and the Proposition is proved.

\qed

\medskip

The following result is a generalization of the previous
proposition.

\begin{prop}\label{prop:chain.com2} %
(1) Any morphism $R^{\ast}:F^{\ast}\rightarrow E^{\ast}$ between
two Banach complexes induces a sequence of linear homomorphisms
$H^{k}R^{\ast}:H^{k}(F^{\ast})\rightarrow H^{k}(E^{\ast})$ from
the  cohomology of $F^*$ to the cohomology of $E^*$.

\medskip

(2) The morphism $R^{\ast}:F^{\ast}\rightarrow E^{\ast}$ induces
a sequence  of bounded operators
$\overline{H}^{k}R^{\ast}:\overline{H}^{k}(F^{\ast})\rightarrow\overline{H}^{k}(E^{\ast})$
from the reduced cohomology of $F^*$ to the reduced cohomology of
$E^*$.

\medskip

(3) If there exists a  homotopy between two morphisms $R^{\ast}$
and $S^{\ast}:F^{\ast}\rightarrow E^{\ast}$, then the
corresponding homomorphisms on the cohomology groups coincide:
$$H^{k}R^*=H^{k}S^{\ast}:H^{k}(F^{\ast})\rightarrow
H^{k}(E^{\ast}).$$

\medskip

(4) If there exists a weak  homotopy between two morphisms
$R^{\ast}$ and $S^{\ast}:F^{\ast}\rightarrow E^{\ast}$, then the
corresponding morphisms on the reduced cohomology groups coincide:
$$\overline{H}^{k}R^{\ast}=
\overline{H}^{k}S^{\ast}:\overline{H}^{k}(F^{\ast})\rightarrow\overline{H}^{k}(E^*).$$
\end{prop}

\medskip

\textbf{Proof} \ (1) Because $dR^{*}=R^{*}d$, the image
$R^{*}([\omega])$ of any cohomology class  $[\omega]$ of the
complex $F^{*}$ is a well defined cohomology  class of the complex
$E^{*}$. \\

\medskip

\ (2) Using the continuity of $R^*$ and $dR^{*}=R^{*}d$, we see
that closure of the image $R^{*}([\omega])$ of a reduced cohomology class of
$F^{*}$ is a well defined reduced cohomology class of $E^{*}$. By the
boundedness of  $R^{k}$, the operators
$\overline{H}^{k}R^{\ast}:\overline{H}^{k}(F^{\ast})\rightarrow\overline{H}^{k}(E^{\ast})$
is also bounded.

\medskip

(3) The condition $S^{k}-R^{k}=d\circ A^{k}+A^{k+1}\circ d$ \
implies that for any $\xi\in Z^{k}(F^{*})$ we have
$\left(S^{k}\xi-R^{k}\xi\right)= d(A^{k}\xi)\in B^{k}(E^{*})$.

\medskip

(4) The condition $\lim_{j\rightarrow\infty}\left\Vert (d\circ
A_{j}^{k}+A_{j}^{k+1}\circ d)x-(S^{k}-R^{k})x\right\Vert =0$ for
any $x\in F^{k}$ implies that for any $\xi\in Z^{k}(F^{*})$ we
have $$\lim_{j\rightarrow\infty}\left\Vert
S^{k}\xi-R^{k}\xi-d(A_{j}^{k}\xi)\right\Vert =0.$$ \qed

\medskip

A special case of the previous Proposition is given in the following
definitions:

\begin{defn}
 \textbf{a)} A Banach complex $F^{\ast}=\{ F^{k},d\}$
 is \emph{acyclic} if there exists a family of
 bounded operators $A^{k}: F^{k}\rightarrow F^{k-1}$
 such that
\[
Id=d\circ A^{k}+A^{k+1}\circ d.
\]
 \textbf{b)} The Banach complex $F^{\ast}$ is \emph{weakly acyclic}
if for any $k$ there exists a sequence of bounded operators
$A_{j}^{k}:F^{k}\rightarrow F^{k-1}$ such that for any element
$x\in F^{k}$ we have \[ \lim_{j\rightarrow\infty}\left\Vert (d\circ
A_{j}^{k}+A_{j}^{k+1}\circ d)x-x\right\Vert =0.\]
\end{defn}

\bigskip  

In other words, $F^{\ast}$ is (weakly) acyclic if and only if there exists a (weak) homotopy from the
identity $Id : F^{*} \to F^{*}$ to the trivial morphism $0 : F^{*} \to F^{*}$ 
It is thus clear that an acyclic complex has trivial cohomology and a weakly
acyclic complex has trivial reduced cohomology.


\section{$L_{q,p}$-cohomology and Banach complexes}

In this section, we explain how the $L_{q,p}$-cohomology of a
Riemannian manifold $(M,g)$ can be formally seen as the
cohomology of some complex of Banach spaces. Let us start by
introducing the notation
\[
\Omega_{q,p}^{k}(M):=\left\{ \,\omega\in L^{q}(M,\Lambda^{k})\,
\big| \, d\omega\in L^{p}\right\}.
\]
This is a Banach space for the graph norm
\begin{equation}\label{grafnorm}
\left\Vert \omega\right\Vert _{\Omega_{q,p}}:=\left\Vert
\omega\right\Vert _{L^{q}}+\left\Vert d\omega\right\Vert _{L^{p}}.
\end{equation}

By standard arguments of functional analysis  (see e.g. 
\cite{Brezis}) , it can be proved that $\Omega_{q,p}^{k}(M)$ is a
reflexive Banach space for any $1<p,q<\infty$. We will also prove
in section \ref{sec.regularization} that smooth forms are dense in
$\Omega_{q,p}^{k}(M)$ for any $1\leq p,q<\infty$.

\bigskip

To define a Banach complex, we
 choose an arbitrary finite sequence of numbers
\[
\pi=\{p_{0},p_{1},\cdots,p_{n}\}\subset[1,\infty],
\]
and define
\[
\Omega_{\pi}^{k}(M):=\Omega_{p_{k},p_{k+1}}^{k}(M).
\]

Observe that $\Omega_{\pi}^{n}(M)=L^{p_{n}}(M,\Lambda^{n})$ and
$\Omega_{p,p}^{1}(M)$ coincides with the Sobolev space $W^{1,p}(M)$.

Since the exterior differential is a bounded operator \
$d:\Omega_{\pi}^{k-1}\rightarrow\Omega_{\pi}^{k}$, we have
constructed a Banach complex.
\[
0\rightarrow\Omega_{\pi}^{0}\overset{d}{\rightarrow}\cdots\overset{d}{\rightarrow}\,\,
\Omega_{\pi}^{k-1}\overset{d}{\rightarrow}\,\,
\Omega_{\pi}^{k}\overset{d}{\rightarrow}\cdots\overset{d}{\rightarrow}\
\Omega_{\pi}^{n}\rightarrow0\,.
\]

\begin{defn}
\emph{The (reduced) $L_{\pi}$-cohomology of $M$ is the (reduced)
cohomology of the Banach complex} $\{\Omega_{\pi}^{k}(M),d_{k}\}$\emph{.}
\end{defn}

The $L_{\pi}$-cohomology space $H_{\pi}^{k}(M)$ depends only on
$p_{k\text{ }}$and $p_{k-1}$ and we have in fact
\[
 H_{\pi}^{k}(M)=H_{p_{k-1},p_{k}}^{k}(M) \quad \text{ and } \quad
 \overline{H}_{\pi}^{k}(M) = \overline{H}_{p_{k-1},p_{k}}^{k}(M).
\]

Two cases are of special interest:

\begin{enumerate}
\item The $L_{p}$-cohomology, which corresponds to the constant sequence
$\pi=\{ p,p,...,p\}$.
\item The \emph{conformal cohomology}, which corresponds to the sequence
$p_{0}=\infty$, and $p_{k}=\frac{n}{k}$ for $k=1,...,n$. The
cohomology associated to this sequence is a  conformal invariant
of the manifold by Theorem \ref{th:conf}. \\
Let us remark here
that $\left(\frac{1}{p_{k}}-\frac{1}{p_{k-1}}\right)=\frac{1}{n}$.
\end{enumerate}



\section{$L_{q,p}$-cohomology and Sobolev inequality}

We are now in position to give the  interpretation of $L_{q,p}$-cohomology 
in terms of a Sobolev type inequality for
differential forms on a Riemannian manifold $(M,g)$:

\begin{thm} \label{th.sobin1a}
\ $H_{q,p}^{k}(M,g)=0$
 if and only if there exists a constant $C<\infty$ such that for
any closed $p$-integrable differential form $\omega$ of degree
$k$ there exists a differential form $\theta$ of degree $k-1$ such
that $d\theta=\omega$ and
\[
\left\Vert \theta\right\Vert _{L^{q}}\leq C\left\Vert \omega\right\Vert _{L^{p}}.
\]
\end{thm}

This result is a direct consequence of Proposition \ref{abs.sob1}.

\qed

\begin{thm}\label{th.sobin2a}
A) If $T_{q,p}^{k}(M)=0$, \,
  then there exists a constant $C'$ such that for any differential
form $\theta \in \Omega^{k-1}_{q,p}(M)$ of degree $k-1$ there exists a
closed form $\zeta\in Z_{q}^{k-1}(M)$ such that
\begin{equation}\label{inconc.sob2}
\left\Vert \theta-\zeta\right\Vert _{L^{q}} \leq C'\left\Vert
d\theta\right\Vert _{L^{p}}.
\end{equation}

B) Conversely, if $1<q<\infty$, and if there exists a constant
$C'$ such that for any form $\theta \in \Omega^{k-1}_{q,p}(M)$ of
degree $k-1$ there exists  $\zeta\in Z_{q}^{k-1}(M)$ such that
(\ref{inconc.sob2}) holds, then $T_{q,p}^{k}(M)=0$.
\end{thm}
This statement follows immediately from  Proposition
\ref{abs.sob2} and \ref{abs.sob3}.

\qed


\section{Manifolds with finite volume and monotonicity}
\label{sec.mon1}

The $L_{q,p}$-cohomology of a manifold with finite volume has some monotonicity properties. In the next statement, 
the symbol $H_{2}\twoheadrightarrow H_{1}$ (where $H_{1},H_{2}$ are vector spaces) means that $H_{1}$ is a quotient of $H_{2}$.

\medskip

\begin{prop}
If $(M,g)$ has finite volume, $\ 1\leq p \leq \infty $ and $\ 1\leq q_{1}\leq q_{2}\leq \infty $, 
then $\overline{H}_{q_{2},p}^{k}(M)\twoheadrightarrow \overline{H}_{q_{1},p}^{k}(M)$ and 
${H}_{q_{2},p}^{k}(M)\twoheadrightarrow  {H}_{q_{1},p}^{k}(M)$.
\end{prop}

\textbf{Proof} 
Since $1\leq q_{1}\leq q_{2}$ and $M$ has finite volume, we have
$L^{q_{1}}(M,\Lambda^{k})\supset L^{q_{2}}(M,\Lambda^{k})$, hence
$\Omega^{k-1}_{q_{1},p}\supset \Omega^{k-1}_{q_{2},p}$ and thus 
\begin{eqnarray*}
\overline{B}_{q_1,p}^{k}(M) &=& \overline{d\left(
 \Omega^{k-1}_{q_{1},p}\right)}\cap L^{p}(M,\Lambda^{k})
 \\ & \supset &
\overline{d\left( \Omega^{k-1}_{q_{2},p}\right)}\cap L^{p}(M,\Lambda^{k})
\\  &=& \overline{B}_{q_2,p}^{k}(M) .
\end{eqnarray*}
Since $B_{2}\subset B_{1} \subset Z$ implies $Z/B_{1}\twoheadrightarrow   Z/B_{2}$, we have
$$
 \overline{H}_{q_{2},p}^{k}(M) =
 Z^k_p/ \overline{B}_{q_2,p}^{k}(M)
\twoheadrightarrow 
Z^k_p/ \overline{B}_{q_1,p}^{k}(M) =
\overline{H}_{q_{1},p}^{k}(M).
$$
The proof for unreduced cohomology is the same.

\qed

\medskip

We also have some kind of monotonicity with respect to $p$:

\begin{prop}
If $(M,g)$ has finite volume $\ 1\leq p_{2}\leq p_{1}\leq \infty $
\ and $\ 1\leq q_{1}\leq q_{2}\leq \infty $, then
$$
 H_{q_{2},p_2}^{k}(M)=0  \ \Rightarrow  \ H_{q_{1},p_1}^{k}(M) = 0.
$$
\end{prop}

\textbf{Proof}  Since  $M$ has finite volume, $q_{1}\leq q_{2}$
and $p_{2}\leq p_{1}$, we have\footnote{The symbol $\lesssim$
means that the inequality holds up to some constant.} for any
$q_2$-integrable form $\theta$ and  any $p_1$-integrable form
$\omega$

\[
 \left\Vert \theta\right\Vert _{L^{q_1}} \lesssim \left\Vert  \theta\right\Vert _{L^{q_2}}
 \qquad \text{ and } \qquad
  \left\Vert \omega\right\Vert _{L^{p_2}} \lesssim \left\Vert \omega\right\Vert _{L^{p_1}}.
\]

Since $H_{q_{2},p_2}^{k}(M)=0$, we know from Theorem \ref{th.sobin1a} that  for any closed $p_2$-integrable form $\omega$ of degree $k$ there
exists a differential form $\theta$ of degree $k-1$ such that
$d\theta=\omega$ and
$$ 
\left\Vert \theta\right\Vert _{L^{q_2}} \lesssim \left\Vert  \omega \right\Vert
_{L^{p_2}}.
$$
Combining this inequality with two previous inequalities we get
$$ \left\Vert \theta\right\Vert _{L^{q_1}} \lesssim \left\Vert  \omega \right\Vert
_{L^{p_1}}$$ and the result immediately follows from the same
Theorem \ref{th.sobin1a}.

\qed

\medskip

For the torsion, we need to avoid the values $q=1$ and $q=\infty$:

\begin{prop}
If $(M,g)$ has finite volume $\ 1\leq p_{2}\leq p_{1}\leq \infty $
\ and $\ 1< q_{1}\leq q_{2} < \infty $, then
$$
 T_{q_{2},p_2}^{k}(M)=0  \ \Rightarrow  \ T_{q_{1},p_1}^{k}(M) = 0.
$$
\end{prop}

\textbf{Proof}  Again, since  $q_{1}\leq q_{2}$ we have
$\zeta\in Z_{q_2}^{k-1}(M)
 \Rightarrow \zeta\in Z_{q_1}^{k-1}(M)$ and
\[
 \left\Vert \theta-\zeta\right\Vert _{L^{q_1}} \lesssim \left\Vert  \theta-\zeta\right\Vert _{L^{q_2}}
 \qquad \text{ and } \qquad
  \left\Vert d\theta\right\Vert _{L^{p_2}} \lesssim \left\Vert d\theta\right\Vert _{L^{p_1}}.
\]
We may thus argue as in the previous proof using Theorem \ref{th.sobin2a}.

\qed


\section{Almost duality} \label{sec.almostduality}

It  has been proved in \cite{GK4}  that for complete manifolds
the dual space of $\overline{H}_{p}^{k}(M)$ coincides with
$\overline{H}_{p^{'}}^{n-k}(M)$ where
$\frac{1}{p}+\frac{1}{p^{'}}=1$ (there is also a duality result
for non complete manifolds). The duality is based on the pairing
$\int_{M}\alpha\wedge\beta$ where $\alpha\in\Omega_{p}^{k}(M)$
and $\beta\in\Omega_{p^{'}}^{k}(M)$.

\medskip

For $L_{q,p}$-cohomology we have no convenient description of
dual spaces, but the notion of \emph{almost duality} which we now
introduce is sufficient for many calculations.

\medskip

We start with a rather elementary result about the non vanishing of
$L_{q,p}$-cohomology:

\begin{lem}
Let $(M,g)$ be an arbitrary Riemannian manifold of dimension $n$.
Let \ $\alpha \in Z_{p}^{k}(M)$.  If there exists $\gamma \in
C_{\text{c}}^{\infty }(M,\Lambda ^{n-k})$ such that $d\gamma =0$\
and \ $\int_{M}\alpha \wedge \gamma \neq 0$, then $[\alpha ]\neq
0$ in $\overline{H}_{q,p}^{k}(M)$ for any $1\leq q\leq \infty$.
\end{lem}

\textbf{Proof} \ Suppose that $\alpha \in
\overline{B}_{q,p}^{k}(M)$. Then $\displaystyle\alpha
=\lim_{j\rightarrow \infty }d\beta _{j}$ (where the limit is in
$L^{p}$-topology) for some $\beta _{j}\in L^{q}(M,\Lambda ^{k-1})$
with $d\beta _{j}\in L^{p}(M,\Lambda ^{k})$. We then have for any
closed form with compact support $\gamma \in C_{\text{c}}^{\infty
}(M,\Lambda ^{n-k})$
\[
\int_{M}\gamma \wedge \alpha =\lim_{j\rightarrow \infty
}\int_{M}\gamma \wedge d\beta _{j}=\lim_{j\rightarrow \infty
}(-1)^{n-k+1}\int_{M}d\gamma \wedge \beta _{j}=0
\]
in contradiction to the assumption. \qed

\bigskip

There are several generalizations of this result :

\begin{prop} \label{pr.nv1}
Let $(M,g)$ be an arbitrary Riemannian manifold of dimension $n$.
Let \ $\alpha \in Z_{p}^{k}(M)$. Then
\newline
\textbf{A)} If there exists a sequence $\ \{\gamma _{i}\}\subset
C_{\text{c}}^{\infty }(M,\Lambda ^{n-k})$ such that \

\begin{enumerate}[i)]
 \item   $\ds \liminf_{i\rightarrow \infty }  \int_{M}\alpha
\wedge \gamma_{i}>0$;
\item  $\ds \lim_{i\rightarrow \infty }\left\| d\gamma _{i}\right\|_{q'}=0$ \
where $q'=\frac{q}{q-1}$.
\end{enumerate}
Then $[\alpha ]\neq 0$ in $H_{q,p}^{k}(M)$.
\newline

\medskip

\textbf{B)} If there exists a sequence $ \{\gamma _{i}\}\subset
C_{\text{c}}^{\infty }(M,\Lambda ^{n-k})$ satisfying the
conditions (i) and (ii) above and
\begin{enumerate}[i)]\setcounter{enumi}{2}
\item $\left\| \gamma _{i}\right\|_{p'}$ is a
bounded sequence for  $p'=\frac{p}{p-1}$.
\end{enumerate}
Then $[\alpha ]\neq 0$ in $\overline{H}_{q,p}^{k}(M)$.
\end{prop}

\bigskip

\textbf{Proof} \ \textbf{A)}  Suppose that $\alpha =d\beta $ for
some $\beta \in L^{q}(M,\Lambda ^{k-1})$, then by H\"{o}lder
inequality we have for any $\gamma \in C_{\text{c}}^{\infty
}(M,\Lambda ^{n-k})$
\[
\left| \int_{M}\alpha \wedge \gamma \right| =\left| \int_{M}d\beta
\wedge \gamma \right| =\left| \int_{M}\beta \wedge d\gamma \right|
\leq \left\| \beta \right\|_q \cdot \left\| d\gamma \right\|_{q'}.
\]
It follows that for any sequence \ $\{\gamma _{i}\}\subset
C_{\text{c}}^{\infty }(M,\Lambda ^{n-k})$ such that
$\lim_{i\rightarrow \infty }\left\| d\gamma \right\|_{q'}=0$, we
have $\ds \lim_{i\rightarrow \infty }\left| \int_{M}\alpha \wedge
\gamma \right| \leq \lim_{i\rightarrow \infty }\left\| \beta
\right\|_q\cdot \left\| d\gamma _{i}\right\| _{L^{q^{\prime
}}(M)}=0$.

\medskip

\textbf{B)} Suppose that $\alpha \in \overline{B}_{q,p}^{k}(M)$.
Then $\displaystyle\alpha =\lim_{j\rightarrow \infty }d\beta _{j}$
for  $\beta _{j}\in L^{q}(M,\Lambda ^{k-1})$ with $d\beta _{j}\in
L^{p}(M,\Lambda ^{k})$. We have for any $i,j$
\[
\int_{M}\gamma _{i}\wedge \alpha =\int_{M}\gamma _{i}\wedge d\beta
_{j}+\int_{M}\gamma _{i}\wedge (\alpha -d\beta _{j})\  .
\]
For each $j\in \Bbb{N}$, we can find $i=i(j)$ large enough so that
\ $\Vert d\gamma _{i(j)}\Vert _{{q}'}\,\Vert \beta _{j}\Vert
_{q}\leq 1/j$, we thus have
\[
\quad \left| \int_{M}\gamma _{i(j)}\wedge d\beta _{j}\right| \quad
\leq \quad \left| \int_{M}d\gamma _{i(j)}\wedge \beta _{j}\right|
\quad \leq \Vert d\gamma _{i(j)}\Vert _{{q}'}\,\Vert \beta
_{j}\Vert _{q}\leq \frac{1}{j}\,.
\]
On the other hand
\[
\quad \lim_{j\rightarrow \infty }\left| \int_{M}\gamma
_{i(j)}\wedge (\alpha -d\beta _{j})\right| \quad \leq
\lim_{j\rightarrow \infty }\Vert \gamma _{i(j)}\Vert _{{p}^{\prime
}}\,\Vert (\alpha -d\beta _{j})\Vert _{p}=0\,
\]
since \ $\Vert \gamma _{i(j)}\Vert _{{p}'}\,$\ is a bounded
sequence and $\ \Vert (\alpha -d\beta _{j})\Vert _{p}\rightarrow
0$. It follows that $\int_{M}\gamma _{i(j)}\wedge \alpha
\rightarrow 0$ in contradiction to the hypothesis.

\qed

\subsection{The case of complete manifolds}

If $M$ is a complete manifold, we don't need to assume that the form
$\gamma$ from the previous discussion has compact support.

\begin{prop}
\label{pr.complnv1} Assume that $M$ is complete. Let  $\alpha\in
Z_{p}^{k}(M)$, and assume that there exists a smooth closed
$(n-k)$-form $\gamma$ such that $\gamma\in Z_{q'}^{n-k}(M)$,\,
for\, $q'=\frac{q}{q-1},\, \gamma\wedge\alpha\in L^{1}(M)$ and
\[
\int_{M}\gamma\wedge\alpha\neq0,
\]
then $\alpha\,\notin\, B_{q,p}^{k}(M)$. In particular,
$H_{q,p}^{k}(M)\neq\emptyset$.
\end{prop}
This proposition has also version for reduced
$L_{q,p}$-cohomology:

\begin{prop}
\label{pr.complnv2} Assume that $M$ is complete. Let $\alpha\in
Z_{p}^{k}(M)$, and assume that there exists a smooth closed
$(n-k)$-form $\gamma\in Z_{p'}^{n-k}(M)\cap Z_{q'}^{n-k}(M)$,
where $p'=\frac{p}{p-1}$ and $q'=\frac{q}{q-1}$, such that
\[
\int_{M}\gamma\wedge\alpha\neq0,
\]
then $\alpha\,\notin\,\overline{B}_{q,p}^{k}(M)$ where
$q'=\frac{q}{q-1}$. In particular,
$\overline{H}_{q,p}^{k}(M)\neq\emptyset$.
\end{prop}
The proofs are based on the following integration by part lemma:

\begin{lem}
\label{lem.stokes} Assume that $M$ is complete. Let $\beta\in
L^{q}(M,\Lambda^{k-1})$ be such that $d\beta\in L^{p}(M,\Lambda^{k})$,
and $\gamma\in L^{p^{'}}(M,\Lambda^{n-k})$ be such that $d\gamma\in
L^{q^{'}}(M,\Lambda^{n-k+1})$ where
$\frac{1}{p}+\frac{1}{p'}=\frac{1}{q}+\frac{1}{q'}=1$.

If $\gamma$ is smooth and $\gamma\wedge d\beta\in L^{1}(M)$, then
\begin{equation}
\int_{M}\gamma\wedge d\beta=(-1)^{n-k+1}\int_{M}d\gamma\wedge\beta,\label{intbypart}\end{equation}
 In particular, if $\gamma\in L_{p^{'}}^{n-k}(M)\cap L_{q^{'}}^{n-k+1}(M)$,
then the above conclusion holds.
\end{lem}
\textbf{Proof}
The integrability of $d\gamma\wedge\beta$ and $\gamma\wedge d\beta$
is a direct consequence of H\"{o}lder's inequality.

By  H\"{o}lder's inequality, the forms $d\gamma\wedge\beta$ and
$\gamma\wedge d\beta$ both belong to $L^{1}(M)$.

If $\gamma$ is a smooth form with compact support, then the equation
(\ref{intbypart}) follows from the definition of the weak exterior
differential (of $\beta$).

If the support of $\gamma$ is not compact, we set $\gamma_{i}:=\psi_{i}\gamma$
where $\{\psi_{i}\}$ is a sequence of smooth functions with compact
support such that $\psi_{i}(x)\rightarrow1$ uniformly on every compact
subset, $0\leq\psi_{i}(x)\leq1$ and $|d\psi_{i}|_{x}\leq1$ for all
$x\in M$ (such a sequence exists on any complete manifold).

The formula (\ref{intbypart}) holds for each $\gamma_{i}$ (since
these forms have compact support).

Using $|d\psi_{i}|_{x}\leq1$, we have the estimate
\[
|\gamma_{i}\wedge d\beta+(-1)^{n-k}d\gamma_{i}\wedge\beta|\leq|d\gamma\wedge\beta|
+|\gamma\wedge d\beta|+|\gamma\wedge\beta|\in L^{1}(M).
\]
By Lebesgue's dominated convergence theorem, we thus have
\[
\int_{M}\,\left(\gamma\wedge d\beta+(-1)^{n-k}d\gamma\wedge\beta\right)
=\lim_{i\rightarrow\infty}\int_{M}\,\left(\gamma_{i}\wedge
d\beta+(-1)^{n-k}d\gamma_{i}\wedge\beta\right)=0\,.
\]

\qed

\textbf{Proof of Proposition \ref{pr.complnv1}}  Suppose that
$\alpha\in B_{q,p}^{k}(M)$. Then $\alpha=d\beta$ for some
$\beta\in L^{q}(M,\Lambda^{k-1})$. By the previous lemma, we have
\[
\int_{M}\gamma\wedge\alpha=\int_{M}\gamma\wedge
d\beta=(-1)^{n-k+1}\int_{M}d\gamma\wedge\beta=0
\]
(since $\gamma$ is closed) in contradiction to the assumption.
\qed

\bigskip

\textbf{Proof of Proposition \ref{pr.complnv2}}  Suppose that
$\alpha\in\overline{B}_{q,p}^{k}(M)$. Then $\ds
\alpha=\lim_{j\rightarrow\infty}d\beta_{j}$ (where the limit is
in $L^{p}$-topology) for some $\beta_{j}\in
L^{q}(M,\Lambda^{k-1})$ with $d\beta_{j}\in L^{p}(M,\Lambda^{k})$.
Since $d\gamma = 0$, we have
\[
\int_{M}\gamma\wedge\alpha=\lim_{j\rightarrow\infty}\int_{M}\gamma\wedge
d\beta_{j}=\lim_{j\rightarrow\infty}(-1)^{n-k+1}\int_{M}d\gamma\wedge\beta_{j}=0,
\]
which contradicts our hypothesis. \qed


\section{The $L_{q,p}$-cohomology of the line} 
 
 In the following three sections, we compute the $L_{q,p}$-cohomology of the line, the hyperbolic plane and the ball.
 We will see in particular that the only case where  ${H}^1_{q,p}(\mathbb{R})$ vanishes is 
 when  $q= \infty$, $p=1$ :

\begin{prop}
\  $H_{\infty,1}^{1}(\mathbb{R})=0$.
\end{prop}

\textbf{Proof}
If $\omega=a(x)dx$ belongs to $L^{1}(\mathbb{R})$, then $f(x):=\int_{-\infty}^{x}a(s)ds$
belongs to $L^{\infty}(\mathbb{R})$, hence $H_{1,\infty}^{1}(\mathbb{R})=0$.

\qed

\begin{prop} \label{prop.torsR}
 $T_{q,p}^{1}(\mathbb{R})\neq 0$   for any  $1\leq p, q \leq \infty$
 with the  only exception of $q= \infty$, $p=1$.
\end{prop}

\medskip

\textbf{Proof} 
Assume first that $q<\infty$. We know from Theorem \ref{th.sobin2a} that if we had
 $T_{q,p}^{1}(\mathbb{R})=0$, then  there would exist a
Sobolev inequality for functions on the real line $\mathbb{R}$:
\begin{equation}\label{inconc.sobR}
\inf_{z \in \r} \left( \int_{-\infty}^{\infty}|f(x)-z|^q dx\right)^{1/q} \leq C \cdot 
  \left( \int_{-\infty}^{\infty}|f'(x)|^p dx\right)^{1/p}
\end{equation}
 for some constant $C<\infty$.

To see that no such inequality is possible, consider a family of
smooth functions with compact support  $f_a: \mathbb{R} \to \mathbb{R}$  such that $f(x)=1$ if $x\in [1,a]$ and $f_a(x)=0$ if
 $x\not\in [0,a+1]$. We may also assume that $\|f_a'\|_{L^{\infty}}\leq 2$. Assume now that the inequality (\ref{inconc.sobR})
holds. Then the constant $z$ must be zero and we have
$$\int_{-\infty}^{\infty}|f_a(x)|^q dx \geq a-1
\qquad  \text{and} \qquad
  \int_{-\infty}^{\infty}|f_a'(x)|^p dx\leq 2^{1+p},$$
hence
$$C \geq 2^{-1-\frac{1}{p}}(a-1)^\frac{1}{q}$$
for all $a>0$ and we conclude that $C=\infty$.

\medskip

Assume now that  $q =\infty$ and $p>1$. Again, if  we had  $T_{\infty,p}^{1}(\mathbb{R})=0$,  there would exist  $C<\infty$
such that for any  $f\in L^p(\mathbb{R})$:
\begin{equation}\label{inconc.sobRR}
\inf_{z \in \r} \|f(x)-z\|_{\infty} \leq C \cdot \|f'(x)\|_{L^p(\mathbb{R})}.
\end{equation}
Let us consider the functions $g_{k}(x):= e^{-\pi k x^2}$ and $f(x):=\int_{-\infty}^{x}g(u)du$.

We have $0\leq f(x) < \sup f = \int_{-\infty}^{\infty}g(u)du = \frac{1}{\sqrt{k}}$, hence 
$\inf_{z \in \r} \|f(x)-z\|_{\infty} = \frac{1}{2\sqrt{k}}$. On the other hand 
$ \|f'(x)\|_{L^p(\mathbb{R})}= (kp)^{-1/2p}$, hence the constant in (\ref{inconc.sobRR}) satisfies
$$
  \frac{1}{2} k^{-1/2} \leq  C \cdot (kp)^{-1/2p}
$$
for all $k>0$, i.e.  $C=\infty$ since $p>1$.

\medskip

Finally, we have $T_{\infty,1}^{1}(\mathbb{R})=0$ since $H_{\infty,1}^{1}(\mathbb{R})=0$.

\qed

\medskip

Let us turn to the reduced cohomology:

\begin{prop} 
 $\overline{H}^1_{q,p}(\mathbb{R}) \neq 0$ if and only if $p=1$ and  $1\leq q < \infty$.
\end{prop}

\smallskip

\textbf{Proof}
For $p=1,q=\infty$, we have  
$\overline{H}_{\infty,1}^{1}(\mathbb{R}) = H_{\infty,1}^{1}(\mathbb{R})=0$.

\medskip

Assume $1\leq q \leq  \infty$ and  $1 < p \leq \infty$ and 
let $\omega=a(x)dx\in L^{p}(\mathbb{R})$. For each $m\in\mathbb{N}$,
we set $\omega_{m}:=\chi_{[-m,m]}\omega=(\chi_{[-m,m]}(x)a(x))dx$. 
Let us choose a continuous function $\lambda_{m}(x)$ with compact
support in $[0,\infty)$ such that $\int_{\mathbb{R}}\lambda_{m}(x)dx=\int_{-m}^{m}a(x)dx$
and $\left\Vert \lambda_{m}\right\Vert _{L^{p}(\mathbb{R})}<\frac{1}{m}$.

Let $b_{m}(x):= 
\int_{-\infty}^{x}\left(\chi_{[-m,m]}(t)a(t)-\lambda_{m}(t)\right)dt$,
then $b_{m}\in L^q(\mathbb{R})$  (in fact $b_{m}$ has compact support) and $\left\Vert
db_{m}-\omega\right\Vert _{L^{p}(\mathbb{R})}\leq \left\Vert
a\right\Vert _{L^{p}(\mathbb{R}\setminus[-m,m])}+\left\Vert
\lambda_{m}\right\Vert _{L^{p}(\mathbb{R})}\rightarrow0$ as
$m\rightarrow\infty$. This shows that  $\overline{H}^1_{q,p}(\mathbb{R}) = 0$.

\medskip

Assume now that  $p=1$ and  $1\leq q < \infty$ and let $\omega=a(x)dx$ be a 1-form on $\mathbb{R}$ such that  
$\int_{\mathbb{R}}f\omega =1$
and $a(x)$ is smooth with compact support (say  $\mathrm{supp }(a) \subset [1,2]$). Let $f_{j}: \mathbb{R} \to \mathbb{R}$ 
be a sequence of
smooth functions with compact support such that $f_{j}= 1$ on $[1,2]$, $\|f_{j}\|_{L^{\infty}}=1$ and  
$\|f'_{j}\|_{L^{q'}}\leq \frac{1}{j}$
where $q'=q/(q-1)$.

Using Proposition  \ref{pr.nv1}, we see that   $[\omega] \neq 0  \in \overline{H}^1_{q,1}(\mathbb{R})$, because 
 $\omega\in L^{1}(\mathbb{R})$  and the sequence  $\{ f_{j}\} \subset C^{\infty}_{c}(\mathbb{R})$ satisfies the three conditions of
that Proposition.

\qed

\bigskip

\textbf{Remarks}  1.)
In degree 0, the $L_{q,p}$-cohomology is controlled by the volume:
$\overline{H}_{q,p}^{0}(\mathbb{R})= H_{q,p}^{0}(\mathbb{R})=0$ if and only if $p< \infty$ and
$\overline{H}_{q,\infty}^{0}(\mathbb{R})= H_{q,\infty}^{\infty}(\mathbb{R})=\mathbb{R}$.

\medskip

2.) All the results of this section also hold for the half-line $\mathbb{R}_{+}$.

\medskip


\section{The cohomology of the hyperbolic plane}

We treat in this section the case of the hyperbolic plane.

\smallskip

Recall that the hyperbolic plane is the Riemannian manifold
$\mathbb{H}^{2}=\{(u,v)\in\mathbb{R}^{2}:v>0\}$ with the metric
$ds^{2}=v^{-2}(du^{2}+dv^{2})$.

\begin{thm}\label{coho.h2}
For any $q,p\in(1,\infty)$ we have 
\[ 
 \dim (\bar{H}_{q,p}^{1}(\mathbb{H}^{2}))=\infty\,.
 \]

\end{thm}
It will be convenient to introduce new coordinates (the so called
``horocyclic coordinates'') $y:=u$, $z:=-\log(v)$, so that
$\mathbb{H}^{2}=\{(y,z)\in\mathbb{R}^{2}\}$ with ${\displaystyle
ds^{2}=e^{2z}dy^{2}+dz^{2}}$.

\begin{lem}\label{lem.fg} 
There exist two smooth functions $f$ and $g$ on $\mathbb{H}^{2}$ 
such that
\begin{enumerate}[1.)]
  \item $f$ and $g$ are non negative;
  \item $f(y,z)=g(y,z)=0$ if $z\le0$ or $|y|\ge1$;
  \item $df$ and $dg\in L^{r}(\mathbb{H}^{2},\Lambda^{1})$ for any $1<r\le\infty$;
  \item the support of $df\wedge dg$ is contained in $\{(y,z):|y|\le1\,,\,0\le z\le1\}$;
  \item $df\wedge dg\ge0$;
  \item ${\displaystyle \int\!\!\int_{\mathbb{H}^{2}}df\wedge dg=1}$;
  \item ${\displaystyle \frac{\partial f}{\partial y}}$ and
 ${\displaystyle \frac{\partial g}{\partial y}\in L^{\infty}(\mathbb{H}^{2})}$;
  \item ${\displaystyle \frac{\partial f}{\partial z}}$ and
${\displaystyle \frac{\partial g}{\partial z}}$ have compact support.
\end{enumerate}
\end{lem}

\medskip

\textbf{Remark} The forms $df$ and $dg$ cannot have compact support,
otherwise, by Stokes theorem, we would have $\int_{\mathbb{H}^{2}}df\wedge dg=0$.

\medskip

\textbf{Proof}
Choose smooth functions $h_{1}$, $h_{2}$, and $k:\mathbb{R}\rightarrow\mathbb{R}$
with the following properties:

\  1) $h_{1},h_{2}$ and $k$ are $\geq0$;

\  2) $h_{i}(y)=0$ if $|y|\geq1$;

\  3) $h_{1}^{\prime}(y)h_{2}(y)\geq0$ and $h_{1}(y)h_{2}^{\prime}(y)\leq0$
for all $y\in \mathbb{R}$;

\  4) the function $(h_{1}^{\prime}(y)h_{2}(y)-h_{1}(y)h_{2}^{\prime}(y))$
has non empty support;

\  5) $k^{\prime}(z)\geq0$ for all $z$;

\  6) $k(z)=1$ if $z\geq1$ and $k(z)=0$ if $z\leq0$.

\medskip

We set $f(y,z):=h_{1}(y)k(z)$ and $g(y,z):=h_{2}(y)k(z)$. Properties
(1) and (2) of the lemma are then clear. We prove (3) (i.e. that $df\in L^{r}$
for any $1<r\leq\infty$).

Indeed,
\[
df=h_{1}(y)k^{\prime}(z)dz+k(z)h_{1}^{\prime}(y)dy\,.
\]
The first term $h_{1}(y)k^{\prime}(z)dz$ has compact support, and
the second term $k(z)h_{1}^{\prime}(y)dy$ has its support in the
infinite rectangle $Q=\{|y|\leq1\, z\geq0\}$.

Choose $D<\infty$ such that $|k(z)h_{1}^{\prime}(y)|\leq D$ on $\Omega$.
We have \[
{\displaystyle |k(z)h_{1}^{\prime}(y)dy|\leq D\,|dy|=D\, e^{-z}\,,}\]
 thus, since the element of area of $\mathbb{H}^{2}$ is $dA=e^{z}dydz$,
we have \[
\int_{\mathbb{H}^{2}}|k(z)h_{1}^{\prime}(y)dy|^{r}dA\leq D^{r}\int_{Q}e^{-rz}\, e^{z}dydz\leq2CD^{r}\int_{0}^{\infty}e^{(1-r)z}\, dz<\infty\,,\]

from which one gets $df\in L^{r}$.

Now observe that
\[
df\wedge dg=
\left((k(z)k^{\prime}(z))(h_{1}^{\prime}(y)h_{2}(y)-h_{1}(y)h_{2}^{\prime}(y)\right)\, dy\wedge dz\,,
\]
hence the properties (4) and (5) follow from the construction of
$h_{1},h_{2}$ and $k$.

Property (6) is only a normalization. It can be achieved by multiplying
$f$ (or $g$) by a suitable constant.

Properties (7) and (8) are easy to check.

\qed

\medskip

\textbf{Proof of Theorem \ref{coho.h2}} Define the $1$-forms
$\alpha=df$ and $\gamma=dg$ on $\mathbb{H}^{2}$ (where $f$ and
$g$ are as in Lemma \ref{lem.fg}). It is clear that
$d\alpha=d\gamma=0$. We also know that $\alpha\in L^{p}$ for any
$1<p<\infty$ and that $\gamma$ is smooth and $\gamma\in
L^{p'}\cap L^{q'}$ for all $1<p',q'<\infty$.

Since ${\displaystyle
\int_{\mathbb{H}^{2}}\alpha\wedge\gamma\ne0}$, we see by
proposition \ref{pr.complnv2} that
$\alpha\not\in\overline{B}_{q,p}^{1}(\mathbb{H}^{2})$.

Now using the isometry group of $\mathbb{H}^{2}$, we produce an
infinite family of linearly independent classes in
$\overline{H}_{q,p}^{1}(\mathbb{H}^{2})$. \qed

\section{The cohomology of the ball}  \label{sec.Bn}

Since the unit ball $\mathbb{B}^{n}\subset\mathbb{R}^{n}$ has finite
volume, we have for all $1\leq p,q\leq\infty$ $H_{q,p}^{0}(\mathbb{B}^{n})=
\overline{H}_{q,p}^{0}(\mathbb{B}^{n})=\mathbb{R}$.

\medskip

In higher degree, the vanishing of the De Rham cohomology of
$\mathbb{B}^{n}$ is traditionally called the Poincar\'{e} Lemma;
it is proved by explicitly constructing a primitive to any closed
form. To prove the vanishing of the $L_{q,p}-$cohomology of the
ball, we need to control the $L^{q}$ norm of the primitive of a
closed $L^{p}$-norm. For the case $p=q$, this was done by
Gol'dshtein, Kuz'minov and Shvedov in \cite[Lemma 3.2]{GK1} and
for more general $q$ by Iwaniec and Lutoborski in \cite{IL93}.
They proved the following

\begin{thm}\label{iwan}
For any bounded convex domain \ $U\subset\mathbb{R}^{n}$ and any
$k=1,2,...,n$, there exists an operator
\[
T=T_{U}:L_{loc}^{1}(U,\Lambda^{k})\rightarrow
L_{loc}^{1}(U,\Lambda^{k-1})
\]
with the following properties:
\begin{enumerate}[a.)]
\item  $T(d\theta)+dT\theta=\theta$ (in the sense of currents);
\item $\ds
\left|T\theta(x)\right|\leq C\int_{U}\frac{|\theta(y)|}{\;\:|y-x|^{n-1}}dy$.
\end{enumerate}
\end{thm}

\qed

\begin{cor}\label{ciwana}
The operator $T$ maps $L^{p}(U,\Lambda^{k})$ continuously to
$L^{q}(U,\Lambda^{k-1})$ in the following cases: \\
either 
 \begin{enumerate}
\item[i)]  $1\leq p,q \leq \infty$ and $\frac{1}{p}-\frac{1}{q} <
\frac{1}{n}$,
\end{enumerate}
or
\begin{enumerate}
\item[ii)] $1 < p,q \leq \infty$ and $\frac{1}{p}-\frac{1}{q} \leq
\frac{1}{n}$.
\end{enumerate}
\end{cor}

\medskip

\textbf {Remark} Note that condition (i) is equivalent to $p \geq
n$ or $p<n$ and $q<\frac{np}{n-p}$ and condition (ii) is relevant
to conformal cohomology
$\frac{1}{p_k}-\frac{1}{p_{k-1}}=\frac{1}{n}$.

\medskip

\textbf{Proof}  Assume first that $\frac{1}{p}-\frac{1}{q} < \frac{1}{n}$ and recall the  
Young inequality for convolution (see \cite[Prop. 8.9]{folland}), which says that if
$1\leq r,s,t \leq \infty$ satisfy
$\frac{1}{r}+\frac{1}{s} = 1+\frac{1}{t}$, then $\|f*g \|_{L^t} \leq \|f \|_{L^r}\|g \|_{L^s}$.
Applying this inequality to $f=|\theta|$ and $g=|x|^{1-n}$ with $r=p$, $t=q$ and $s=\frac{pq}{p+pq-q}$, and observing that
$$
 \frac{1}{p}-\frac{1}{q} < \frac{1}{n} \quad \Leftrightarrow \quad
 s(1-n)>-n \quad \Leftrightarrow \quad
 \|g \|_{L^s(U)} < \infty ,
$$
we conclude from previous proposition that  $T :L^{p}(U,\Lambda^{k}) \to L^{q}(U,\Lambda^{k-1})$ is bounded with norm at most
$ \||x|^{1-n} \|_{L^s(U)}$.

\smallskip

If $p>1$ and $\frac{1}{p}-\frac{1}{q} = \frac{1}{n}$, then the conclusion also holds by  the 
Hardy-Litlewood-Sobolev inequality (see \cite[p. 119]{stein}).

\qed

\medskip

\begin{cor}\label{ciwan}
 The operator
 $T: \Omega_{p,r}^{k}(U)\rightarrow\Omega_{q,p}^{k-1}(U)$
is bounded and for any $\omega\in\Omega_{p,r}^{k}(U)$ we have
$Td\omega+dT\omega=\omega$ provided either

i)  $1\leq p,q,r \leq \infty$ such that $\frac{1}{p}-\frac{1}{q} <
\frac{1}{n}$ and $\frac{1}{r}-\frac{1}{p} < \frac{1}{n}$, \

or

ii) $1 < p,q,r \leq \infty$ such that $\frac{1}{p}-\frac{1}{q}
\leq \frac{1}{n}$ and $\frac{1}{r}-\frac{1}{p} \leq \frac{1}{n}$.
\end{cor}

\textbf{Proof} The proof is immediate from the previous Theorem and Corollary.

\qed

\medskip

The Corollary \ref{ciwana} implies the following Poincar\'{e} Lemma :

\begin{prop}\label{poinc}
Suppose that   $p,q$ satisfy either 

 \begin{enumerate}
\item[i)]  $1\leq p,q \leq \infty$ and $\frac{1}{p}-\frac{1}{q} <
\frac{1}{n}$,
\end{enumerate}
or
\begin{enumerate}
\item[ii)] $1 < p,q \leq \infty$ and $\frac{1}{p}-\frac{1}{q} \leq
\frac{1}{n}$.
\end{enumerate}
Then $H_{q,p}^{k}(\mathbb{B}^{n})=0$ 
for any  $k=1,...,n$.
\end{prop}

\textbf{Proof} Let $\omega$ be an arbitrary element in
 $Z_{p}^{k}(\mathbb{B}^{n})$. By Corollary \ref{ciwana},  
 we  have $T\omega\in L^{q}(\mathbb{B}^{n}, \Lambda^{k+1})$, 
 since $\omega=dT\omega+Td\omega=d(T\omega)$ we conclude 
 that $[\omega]=0\in H_{q,p}^{k}(\mathbb{B}^{n})$ and thus $H_{q,p}^{k}(\mathbb{B}^{n})=0$.
 
 \qed

\smallskip

\medskip

If $p,q > 1$, we have a necessary and sufficient condition :

\medskip

\begin{thm} If \  $1 <  p,q  \leq \infty$ and $k=1,...,n$, then
$H_{q,p}^{k}(\mathbb{B}^{n})=0$ if and only if
$\frac{1}{p}-\frac{1}{q} \leq \frac{1}{n}$.
\end{thm}

\medskip

\textbf{Proof}  We know from the previous Proposition that the condition is sufficient .

\smallskip

To prove that $H_{q,p}^{k}(\mathbb{B}^{n})\neq0$ if $p<n$ and $q>\frac{np}{n-p}$,
we will use Proposition \ref{pr.nv1}. Let us fix a number $\mu$
in the interval $k-\frac{n}{p}<\mu<k-1-\frac{n}{q}$ (which is possible
since $\frac{1}{p}>\frac{1}{q}+\frac{1}{n}$); and choose two forms
$\theta\in C^{\infty}(\mathbb{S}^{n-1},\Lambda^{k-1})$ and
$\varphi\in C^{\infty}(\mathbb{S}^{n-1},\Lambda^{n-k-1})$
such that
\[
\int_{\mathbb{S}^{n-1}}\varphi\wedge d\theta=1.
\]
For any $0<t<1/4$, we choose a smooth function
$h_{t}:\mathbb{R}\rightarrow\mathbb{R}$ such that $h(t,r)=0$ if
$r<t$ or $r>1-t$ and $h(t,r)=\frac{1}{\left|\log2t\right|}$ if
$r<1-2t$ or $r>2t$.

Let us then consider the forms
\begin{eqnarray*}
\alpha & := & d\left(r^{\mu}\theta\right)\\
\gamma_{t} & := & h_{t}(r)r^{-(\mu+1)}dr\wedge\varphi
\end{eqnarray*}

\smallskip

\textbf{Step 1}  The form $\alpha$ belongs to $L^{p}(\mathbb{B}^{n},\Lambda^{k})$.

\smallskip

We will use the same notation $\theta$ and $\varphi$ for a pullback
of corresponding forms from $\mathbb{S}^{n}$ to $\mathbb{B}^{n}\setminus\{0\}$ induced by 
the radial projection in polar coordinates.

We have
\[
\alpha=r^{\mu}\left(d\theta+\mu\frac{1}{r}dr\wedge\theta\right).
\]
Because $\left|\theta\right|\lesssim r^{-(k-1)}$ and
$\left|d\theta\right|\lesssim r^{-k}$ we have $|\alpha|\lesssim
r^{\mu-k}$. Therefore
\[
\int_{\mathbb{B}^{n}}|\alpha|^{p}dx\lesssim
\int_{0}^{1}\left(r^{\mu-k}\right)^{p}r^{n-1}dr<\infty
\]
 because $p(\mu-k)+n-1>p(k-\frac{n}{p}-k)+n-1>-1$.

\smallskip

\textbf{Step 2} The quantity 
$\left|\int_{\mathbb{B}^{n}}\alpha\wedge\gamma_{t}\right|$ 
is bounded below.

\smallskip

We have
$\alpha\wedge\gamma_{t}=h_{t}(r)r^{-1}dr\wedge\varphi\wedge
d\theta$; since $\int_{\mathbb{S}^{n-1}}\varphi\wedge d\theta=1$,
we have by Fubini Theorem

\[
\left|\int_{\mathbb{B}^{n}}\alpha\wedge\gamma_{t}\right|=
\int_{0}^{1}h_{t}(r)r^{-1}dr \geq
\frac{1}{|\log2t|}\int_{2t}^{1-2t}r^{-1}dr\rightarrow 1
\]
 as $t\rightarrow0$. This implies that  $\left|\int_{\mathbb{B}^{n}}\alpha\wedge\gamma_{t}\right|$ is
bounded below for small values of $t$.

\smallskip

\textbf{Step 3} We have $\left\Vert
d\gamma_{t}\right\Vert _{L^{q'}(\mathbb{B}^{n})}\rightarrow0$ as
$t\rightarrow0$:

\smallskip

We have $d\gamma_{t}:=h_{t}(r)r^{-(\mu+1)}dr\wedge\varphi$ with
$0\leq h_{t}\leq\frac{1}{|\log2t|}$. Since
$\left|dr\wedge\varphi\right|\lesssim r^{-n+k}$,  we have

\[
|d\gamma_{t}|\lesssim \frac{r^{-\mu-1+k-n}}{|\log2t|}
\]
 and by Fubini Theorem
 \begin{eqnarray*}
\int_{\mathbb{B}^{n}}|d\gamma_{t}|^{q'}dx
&= &\int_{\mathbb{B}^{n}}|h_{t}(r)r^{-(\mu+1)}dr\wedge\varphi|^{q'}dx
\\
 &\lesssim &
 \left(\frac{1}{|\log2t|}\right)^{q^{\prime}}
\int_{0}^{1}\left(r^{-\mu-1+k-n}\right)^{q^{\prime}}r^{n-1}dr.
\end{eqnarray*}
 Because
\[
q'(-\mu-1+k-n)+n=q'(-\mu-1+k-n(1-\frac{1}{q'}))=q'(-\mu-1+k-\frac{n}{q})>0
\]
we have
\[
\int_{0}^{1}\left(r^{-\mu-1+k-n}\right)^{q^{\prime}}r^{n-1}dr<\infty.
\]

Therefore
\[
\lim_{t\rightarrow0}\int_{\mathbb{B}^{n}}|d\gamma_{t}|^{q'}dx
\lesssim \lim_{t\rightarrow0}
\left(\frac{1}{|\log2t|}\right)^{q^{\prime}}\int_{0}^{1}
\left(r^{-\mu-1+k-n}\right)^{q^{\prime}}r^{n-1}dr=0
\]

Since $\gamma_{t}$ are smooth forms with compact support, Proposition
\ref{pr.nv1} implies that $[\alpha]\neq0$ in $H_{q,p}^{k}(\mathbb{B}^{n})$.

\qed

\begin{cor}
The conformal cohomology of the hyperbolic space $\mathbb{H}^{n}$
vanishes for any degree $k>1$, i.e.
$$
 H_{\frac{n}{k-1}, \frac{n}{k}}^{k}(\mathbb{H}^{n}) = 0.
$$
\end{cor}

\textbf{Proof}  Since the hyperbolic space $\mathbb{H}^{n}$ is conformally equivalent to the ball $\mathbb{B}^{n} \subset \r^n$, 
this result follows at once from the conformal invariance of conformal cohomology and the previous theorem.

\qed

\begin{rem}
Because $H_{q,p}^{1}(\mathbb{H}^{2})\neq 0$ for any $q,p$, the
previous corollary does not hold for $k=1$. 
\end{rem}


\section{Regularization of forms and cohomology classes}
\label{sec.regularization}

In this section we  investigate two different but related
problems. The first one is a density result for smooth forms in
$\Omega_{q,p}^{*}(M)$ and the second one is a result about
representation of the  cohomology  $H_{q,p}^{*}(M)$ by smooth
forms. We will use the de Rham  regularization method
\cite{derham} and its version for $L_{p}$-cohomology \cite{GK3} in
combination with the results of section \ref{sec.Bn}.

\subsection{Regularization operators for differential forms.}

The standard way of smoothing a function in $\r^n$ is by
convolution with a smooth mollifier. This procedure extends to
differential forms and more generally to any tensor. In his book,
De Rham proposes a clever way of localizing this construction and
grafting it on manifolds.

\medskip

Following De Rham, we associate to any vector $v\in\mathbb{R}^{n}$
the map $s_v : \r^n \to \r^n$ defined by
\begin{equation*}
  s_v(x)= \begin{cases}
  h^{-1}(h(x)+v) & \text{if $\|x\| < 1$}, \\
  x &  \text{if $\|x\|\geq 1$}.
  \end{cases}
\end{equation*}
where $h:\mathbb{B}^{n} \to \r^n$ is a radial diffeomorphism such
that
\begin{equation*}
  h(x)= \begin{cases}
  x & \text{if $\|x\| < 1/3$}, \\
  \frac{1}{\|x\|}\exp (\frac{1}{(1-\|x\|^2)})\cdot x &  \text{if $\|x\|\geq 2/3$}.
  \end{cases}
\end{equation*}

\begin{lem}\label{}
The map $v \to s_v$ defines an action of the group $\r^n$ on the
space $\r^n$ satisfying the following properties:
\begin{enumerate}[a.)]
\item  For every $v\in\mathbb{R}^{n}$, the map
$s_v : \r^n \to \r^n$ is a smooth diffeomorphism;
\item  The mapping $s:\mathbb{R}^{n}\times \r^n\rightarrow \r^n$ is
smooth;
\item    $s_{v}$ is the identity outside of $\mathbb{B}^{n}$;
\item  For every $x\in \mathbb{B}^{n}$ the mapping
  $v \mapsto \alpha_{x}(v):=s_{v}(x)$ is a
diffeomorphism of $\mathbb{R}^{n}$ onto $\mathbb{B}^{n}$.
\end{enumerate}
\end{lem}

\textbf{Proof} For the first two assertions, see \cite{derham}.
The assertions (c) and (d) are obvious. \qed

\bigskip

Let us fix an arbitrary bounded convex domain $U$ such that
$\overline{\mathbb{B}}^{n} \subset U \subset \r^n$. We now define
the regularization operator $
 R_{\epsilon} : L^1_{loc}(U,\Lambda^{k}) \to L^1_{loc}(U,\Lambda^{k})
$
by
\[
R_{\varepsilon}\omega
:=\int_{\mathbb{R}^{n}}s_{v}^{*}(\omega)\rho_{\varepsilon}(v)dv
\]
where $\rho_{\varepsilon}(v) = \rho(v/\varepsilon)$ is a standard
mollifier.

\begin{prop}\label{reloc}
The regularization operator defined
above satisfies the following properties :
\begin{enumerate}[1.)]
  \item For any $\omega \in
  L^1_{loc}(U,\Lambda^{k})$, the form $R_{\epsilon}\omega$ is
  smooth in $\mathbb{B}^n$ and \ $R_{\epsilon}\omega=\omega$ in
  $U\setminus \mathbb{B}^n$;
  \item for any $\omega \in \Omega_{q,p}^{k}(U)$, we have
  $dR_{\varepsilon}\omega = R_{\varepsilon}d\omega$.
  \item For any $1\leq p,q<\infty$ and any $\varepsilon > 0$, the
operator
$$
 R_{\varepsilon} : \Omega_{q,p}^{k}(U)\to \Omega_{q,p}^{k}(U)
$$
is bounded and its norm satisfies $\ds
\lim_{\varepsilon\rightarrow 0}\left\Vert
R_{\varepsilon}\right\Vert_{q,p} = 1$;
  \item For any $1\leq p,q  < \infty$ and any
$\omega\in\Omega_{q,p}^{k}(U)$, we have
$$
\lim_{\varepsilon\rightarrow0} \left\Vert
R_{\varepsilon}^{*}\omega-\omega\right\Vert_{p} = 0.
$$
\end{enumerate}
\end{prop}

\medskip

\textbf{Proof} The first two properties are proved in
\cite{derham}. Property (3) follows from  (2) and \cite[Lemma 2]{GK3}  and 
(4) is a standard property of the regularization.

\qed

\subsection{Homotopy operator}

Given a bounded convex domain $U\subset\mathbb{R}^{n}$ containing
the closed unit ball, we introduce the homotopy 
$$
 A_{\epsilon} :=
 (I-R_{\varepsilon})\circ T_U : \ L^1_{loc}(U,\Lambda^{k}) \to
  L^1_{loc}(U,\Lambda^{k-1}),
$$
where $T_U$ is the operator defined in Theorem \ref{iwan}.
\begin{lem} \label{lem.Aishmtpy}
The operator $A_{\varepsilon}$ is a  homotopy between the
Identity and the regularization operator $R_{\varepsilon}$, i.e.
it satisfies
$$(I-R_{\varepsilon})\,\omega=dA_{\varepsilon}\omega+A_{\varepsilon}d\omega.$$
\end{lem}

\medskip

\textbf{Proof} We know from Theorem \ref{iwan} that
$Td\omega+dT\omega=\omega$ for all $\omega \in
L^1_{loc}(U,\Lambda^{k-1})$, hence we have
\begin{eqnarray*}
  dA_{\varepsilon}\omega+A_{\varepsilon}d\omega  &=&
  d(I-R_{\varepsilon})T\omega+(I-R_{\varepsilon})Td\omega
\\ &=&
dT\omega-dR_{\varepsilon}T\omega + Td\omega -
R_{\varepsilon}Td\omega
\\ &=&
 (dT\omega+Td\omega) - R_{\varepsilon}(dT\omega+Td\omega)
\\ &=&
(I-R_{\varepsilon})(Td\omega+dT\omega)\\
&=& (I-R_{\varepsilon})\omega.
\end{eqnarray*}
\qed

\medskip

\begin{prop}
Let $U\subset\mathbb{R}^{n}$ be a bounded convex domain
containing the closed unit ball. Then
$A_{\varepsilon}:\Omega_{p,r}^{k}(U)\rightarrow\Omega_{q,p}^{k-1}(U)$
is a bounded operator for any $k=1,2,...,n$ in the following two cases:

i)  $1\leq p,q,r \leq \infty$ such that $\frac{1}{p}-\frac{1}{q}
< \frac{1}{n}$, \

ii) $1 < p,q \leq \infty$ and $\frac{1}{p}-\frac{1}{q} \leq
\frac{1}{n}$ and $\frac{1}{r}-\frac{1}{p} \leq \frac{1}{n}$.

Furthermore, we have $(I - R_{\varepsilon})\omega =
dA_{\varepsilon}\omega + A_{\varepsilon}d\omega$ for any
$\omega\in\Omega_{p,r}^{k}(U)$ and $A_{\varepsilon}\omega=0$
outside the unit ball.
\end{prop}

\medskip

\textbf{Proof} The first assertion follows from  Proposition
\ref{reloc} and Corollary \ref{ciwan} and the second one is the
previous Lemma. The last assertion follows from the fact that
$R_{\varepsilon} = I$ outside of the unit ball.

\qed

\medskip

\subsection{Globalization}

This regularization operators $R_{\varepsilon}$ and
$A_{\varepsilon}$ can be globalized as follow: given a Riemannian
manifold $(M,g)$, we can find a countable atlas $\{\varphi_{i}:
V_{i}\subset M \to U_i\}_{i\in \mathbb{N}}$ such that $U_i \subset
\r^n$ is a bounded convex domain satisfying
$\overline{\mathbb{B}}^{n} \subset U_i \subset \r^n$ for all $i$
and that $\{B_i\}$ is a covering of $M$, where $B_i :=
\varphi_i^{-1}(B)\subset V_i$. We also assume that $\{V_i\}$ (and
hence $\{B_i\}$) is a locally finite covering of $M$ (we can in
fact assume that any collection of $n+2$ different charts $V_{i}$
has an empty intersection, where $n=\dim M$.)

\medskip

For any $m\in \mathbb{N}$, we define two operators
$$
 R^{(m)}_{\varepsilon} , A^{(m)}_{\varepsilon} :
 L^1_{loc}(M,\Lambda^{m}) \to L^1_{loc}(M,\Lambda^{m})
$$
as follow:
$$
  R^{(m)}_{\varepsilon} :=
  R_{1,\varepsilon}\circ R_{2,\varepsilon}\circ \cdots \circ
  R_{m,\varepsilon},
$$
and
$$
  A^{(m)}_{\varepsilon} :=
  R_{1,\varepsilon}\circ R_{2,\varepsilon}\circ \cdots \circ R_{m-1,\varepsilon}
  \circ A_{m,\varepsilon},
$$
where
$$
 R_{i,\varepsilon} (\theta) := \left(\varphi_i^{-1}\right)^*
 \circ R_{\varepsilon}\circ \varphi_i^* (\theta);
$$
and
$$
 A_{i,\varepsilon} (\theta) := \left(\varphi_i^{-1}\right)^*
 \circ  (R_{i, \varepsilon}-I)\, T_{U_i} \circ \varphi_i^* (\theta).
$$
Here $T_{U_i}$ is the operator defined on the domain $U_i$ in
Theorem \ref{iwan}.

\medskip

Observe that the operator $R_{i,\varepsilon}$ is a priori only
defined on $V_i$, but it acts as the identity on
$V_{i}\setminus\overline{B}_{i}$ and can thus be extended on the
whole of $M$ by declaring that $R_{i,\varepsilon}=id$ on
$M\setminus\overline{B}_{i}$. Likewise, the operator
$A_{i,\varepsilon}$ is a priori only defined on $V_i$, but it is
zero on $V_{i}\setminus\overline{B}_{i}$ (because
$R_{\varepsilon} = I$ outside of the unit ball). Hence
$A_{i,\varepsilon}$ can be extended on the whole of $M$ by
declaring $A_{i,\varepsilon}=0$ on $M\setminus\overline{B}_{i}$.

\medskip

We now define the global regularization operator and the global
homotopy operator as follow:
\begin{equation}\label{def.glob}
   R^M_{\varepsilon}:= \lim_{m\to \infty}R^{(m)}_{\varepsilon},
   \qquad
   A^M_{\varepsilon}:= \sum_{m=1}^{\infty} A^{(m)}_{\varepsilon}.
\end{equation}

\medskip

By construction, the expressions $R^M_{\varepsilon}:= \prod_{i}
R_{i,\varepsilon}$ and $A^M_{\varepsilon}:= \sum_{l}
A^{(k)}_{\varepsilon}$ are really  finite operations in any
compact set and the operators
$R^M_{\varepsilon},A^M_{\varepsilon}$ are thus well defined on
$L^1_{loc}(M,\Lambda^{k})$.

\begin{thm}\label{reglob}
For every Riemannian manifold $M$ there exists a family of
regularization operators $R^M_{\varepsilon}$ and  homotopy
operators $A^M_{\varepsilon}$ such that
\begin{enumerate}[1.)]
   \item For any $\omega \in
  L^1_{loc}(M,\Lambda^{k})$, the form $R^M_{\epsilon}\omega$ is
  smooth in $M$;
  \item For any $\omega \in \Omega_{q,p}^{k}(M)$, we have
  $dR^M_{\varepsilon}\omega = R^M_{\varepsilon}d\omega$;
  \item For any $1\leq p,q<\infty$ and any $\varepsilon > 0$, the
operator $R^M_{\varepsilon} : \Omega_{q,p}^{k}(M)\to
\Omega_{q,p}^{k}(M)$ is bounded and its norm satisfies $\ds
\lim_{\varepsilon\rightarrow 0}\left\Vert
R^M_{\varepsilon}\right\Vert_{q,p} = 1$;
  \item For any $1\leq p,q < \infty$ and any
$\omega\in\Omega_{q,p}^{k}(M)$
 we have
$$
\lim_{\varepsilon\rightarrow 0} \left\Vert
R_{\varepsilon}^{M}\omega-\omega\right\Vert_{p}=0.
$$
\item  The operator   $A_{\varepsilon}:\Omega_{pr}^{k}(M)\rightarrow\Omega_{q,p}^{k-1}(M)$
is bounded for any $k=1,...,n$ in the following cases:
\begin{enumerate}[(i)]
 \item $1\leq p,q,r \leq \infty$ such that $\frac{1}{p}-\frac{1}{q} <
\frac{1}{n}$ and $\frac{1}{r}-\frac{1}{p} < \frac{1}{n}$, 
\item $1 < p,q,r \leq \infty$ such that $\frac{1}{p}-\frac{1}{q}
\leq \frac{1}{n}$ and $\frac{1}{r}-\frac{1}{p} \leq \frac{1}{n}$.
\end{enumerate}
\item  We have the homotopy formula $$\omega - R^M_{\varepsilon}\omega =
dA^M_{\varepsilon}\omega+A^M_{\varepsilon}d\omega. $$
\end{enumerate}
\end{thm}

\textbf{Proof} The first four assertions follow immediately from
Proposition \ref{reloc}.

The fifth assertion follows from  Proposition \ref{reloc} and
Corollary \ref{ciwan}.

To prove the last assertion, observe that by Lemma
\ref{lem.Aishmtpy}, we have $ \omega  - R_{m,\varepsilon}\omega =
dA_{m,\varepsilon}\omega+A_{m,\varepsilon}d\omega $. Multiplying
this expression by $R^{(m-1)}_{\varepsilon}$, we obtain
$$
R^{(m-1)}_{\varepsilon}\omega  - R^{(k)}_{\varepsilon}\omega =
dA^{(k)}_{\varepsilon}\omega+A^{(m)}_{\varepsilon}d\omega,
$$
summing this identities on $m=1,2,...$, we obtain the assertion (6).

\qed

\medskip

\begin{cor}\label{cor.smoothisdense}
For any $q,p \in  [1,\infty)$, the space
$$C^{\infty}\Omega_{q,p}^{k}(M):=C^{\infty}(M)\cap\Omega_{q,p}^{k}(M)$$
of smooth $k$-forms $\theta$\,  in $L^p$ such that $d\theta\in
L^q$ is dense in $\Omega_{q,p}^{k}(M)$.
\end{cor}

\medskip

\textbf{Proof} This result follows immediately from the first
three conditions in Theorem \ref{reglob}. \qed

\subsection{$L_{\pi}$-cohomology and smooth forms}

The previous theorem implies that under suitable assumptions on
$p,q$, the $L_{\pi}$-cohomology of a Riemannian manifold can be
represented by smooth forms.

\medskip

To be more precise, for any sequence $\pi$, we denote by
$$C^{\infty}\Omega_{\pi}^{k}(M):=C^{\infty}(M)\cap\Omega_{\pi}^{k}(M)$$
the subcomplex of smooth forms in $\Omega_{\pi}^{k}(M)$ and by
$$C^{\infty}H^*_{\pi}(M) = H^*(C^{\infty}\Omega_{\pi}^{k}(M))$$
its cohomology.

\medskip

\begin{thm} \label{th.smoothcohomology}
Let $(M,g)$ be a $n$-dimensional
 Riemannian manifold and \\
  $\pi=\{ p_{0},p_{1},\cdots,p_{n}\}\subset (1,\infty)$
a finite sequence of numbers such that \\
$\frac{1}{p_k}-\frac{1}{p_{n-k}} \leq \frac{1}{n}$ for
$k=1,2,..n$. \  Then
$$C^{\infty}H^*_{\pi}(M) = H^*_{\pi}(M).$$
\end{thm}

\textbf{Proof} This result follows immediately from Proposition
\ref{prop:chain.com} and Theorem \ref{reglob}.

\qed

\bigskip

It is perhaps useful to reformulate this theorem without the
language of complexes:

\begin{thm} Let $(M,g)$ be a $n$-dimensional
Riemannian manifold and suppose that $p,q\in (1,\infty)$ satisfy
$\frac{1}{p}-\frac{1}{q} \leq \frac{1}{n}$. Then the cohomology
$H^*_{q,p}(M)$ can be represented by smooth forms.

More precisely, any closed form in $Z^k_p(M)$ is cohomologous to a
smooth form in $L^p(M)$. Furthermore, if two smooth closed forms
$\alpha,\beta \in C^{\infty}(M)\cap Z^k_p(M)$ are cohomologous
modulo $d\Omega_{q,p}^{k-1}(M)$, then they are cohomologous modulo
$dC^{\infty}\Omega_{q,p}^{k-1}(M)$.
\end{thm}

\bigskip

\begin{cor} Let $(M,g)$ be a $n$-dimensional
Riemannian manifold and suppose that $p,q\in (1,\infty)$ satisfy
$\frac{1}{p}-\frac{1}{q} \leq \frac{1}{n}$.
Then any reduced cohomology class can be represented by  a smooth form.
\end{cor}

\textbf{Proof} This is clear from the previous Theorem, since
$\overline{H}^k_{q,p}(M)$ is a quotient of $H^k_{q,p}(M)$. 

\qed

\medskip

\subsection{The case of compact manifolds}

From previous results, we now immediately have:

\medskip

\begin{thm}\label{th.compactLqpcohom}%
Let $(M,g)$ be a compact $n$-dimensional Riemannian manifold and
$\pi=\{ p_{0},p_{1},\cdots,p_{n}\}\subset (1,\infty)$ a finite
sequence of numbers such that $\frac{1}{p_k}-\frac{1}{p_{n-k}}
\leq \frac{1}{n}$ for $k=1,2,..n$. Then
$$H^*_{\pi}(M) = H^*_{\mathrm{DeRham}}(M).$$
In particular $H^*_{\pi}(M)$ is finite dimensional and thus
$T^*_{\pi}(M)=0$.
\end{thm}

\textbf{Proof} Recall that the De Rham cohomology
$H^*_{\mathrm{DeRham}}(M)$ of $M$ is the cohomology of the complex
$(C^{\infty}(M,\Lambda^*),d)$. Any smooth form on a compact
Riemannian manifold clearly belongs to $L^p$ for any $p\in
[0,\infty]$, hence $(C^{\infty}(M,\Lambda^*),d) =
C^{\infty}\Omega_{\pi}^{k}(M)$ and by  Theorem
\ref{th.smoothcohomology}, we have
$$H^*_{\pi}(M)=C^{\infty}H^*_{\pi}(M) = H^*_{\mathrm{DeRham}}(M).$$
It is well known that the De Rham cohomology of a compact manifold
is finite dimensional. Since $\dim T^*_{\pi}(M)\leq \dim
H^*_{\pi}(M) < \infty$, it follows from Lemma \ref{lem.notorsion}
that $T^*_{\pi}(M)=0$.

\qed

\subsection{Proof of Theorems \ref{th.sobin1} and \ref{th.sobin2}}
Let us define the sequence \\ $\pi=\{ p_{0},p_{1},\cdots,p_{n}\}$ by
$p_j=q$ if $j=1,2,..k-1$ and $p_j=p$ if $j=k,...,n$.

By hypothesis, we have $\frac{1}{p}-\frac{1}{q} \leq
\frac{1}{n}$, hence the sequence $\pi$ satisfies
$\frac{1}{p_j}-\frac{1}{p_{j-1}} \leq \frac{1}{n}$ for all $j$.
Hence we know by Theorem \ref{th.compactLqpcohom}  that
$H^k_{q,p}(M) = H^k_{\mathrm{DeRham}}(M)$ and $T^k_{q,p}(M)=0$.

Thus Theorem \ref{th.sobin1} follows from Theorem \ref{th.sobin2a}
and Theorem \ref{th.sobin2} follows from Theorem \ref{th.sobin1a}.

\qed

\section{Relation with a non linear PDE} \label{sec.PDE}

We show in this section that the vanishing of torsion gives
sufficient condition to solving the non linear equation

\begin{equation}\label{smlapl}
 \delta(\left\Vert d\theta\right\Vert ^{p-2}d\theta)=\alpha ,
\end{equation}
where $\delta$ is the operator defined for  $\omega\in  L^{1}_{loc}(M,\Lambda^{k})$ as
$$ \delta \,  \omega= (-1)^{nk+n+1}*d* \, \omega. $$
Recall that for any   $k$-form $\omega$, we have%
\footnote{Here is the proof: Since $\omega$ is a $k$ form, $d*\omega$ is a form of degree $m=n-k+1$ and
$ **d*\omega = (-1)^{m(n-m)}d*\omega = (-1)^{nk+n+1+k}d*\omega$,
therefore
$
(-1)^kd*\omega  = (-1)^{nk+n+1}**d*\omega
= *\delta \omega
$.}.
\begin{equation}\label{eq1*}
 *\delta \omega = (-1)^kd*\omega.
\end{equation}

\medskip

This operator is the  formal adjoint to the exterior differential $d$ in the sense that
\begin{equation}\label{adjointd}
\int_M\langle\omega,d\varphi\rangle \, d\mathrm{vol} =
\int_M\langle \delta\omega,\varphi\rangle \, d\mathrm{vol}
\end{equation}
for any $\varphi\in C_c^{\infty}(M,\Lambda^{k-1})$.

\medskip

Indeed, by definition of the Hodge $*$ operator, we have
$$
 \langle d\varphi , \omega  \rangle \, d\mathrm{vol}
=  \left(  d\varphi \wedge  *\omega \right)
$$
and from the definition of the weak exterior differential, it follows that
\begin{equation*}
\int_M  \langle d\varphi , \omega  \rangle \, d\mathrm{vol} =
\int_M d\varphi \wedge  *\omega
=
(-1)^k\int_M \varphi \wedge  d*\omega .
\end{equation*}
Thus from  (\ref{eq1*}):
\begin{eqnarray*}
\int_M  \langle d\varphi , \omega  \rangle \, d\mathrm{vol}& = &
(-1)^k\int_M \varphi \wedge  d*\omega
\\ &=&
\int_M \varphi \wedge *\delta\omega
\\ &=&
\int_M\langle\varphi , \delta \omega
\rangle \, d\mathrm{vol}.
\end{eqnarray*}

\medskip

Applying (\ref{adjointd}) to $\omega = |d\theta|^{p-2}d\theta$, we obtain the following
\begin{lem}
$\theta\in  L^{1}_{loc}(M,\Lambda^{k})$ is a solution to (\ref{smlapl}) if and only if
\begin{equation}\label{wlapl}
 \int_M  \langle d\varphi , \left\Vert d\theta\right\Vert
^{p-2}d\theta \rangle \, d\mathrm{vol}
=
\int_M\langle\varphi , \alpha \rangle \, d\mathrm{vol}
\end{equation}
for any $\varphi\in C_c^{\infty}(M,\Lambda^{k})$ .
\end{lem}
\qed

The equation (\ref{wlapl}) is just the weak form of  (\ref{smlapl}).

\bigskip

\textbf{Remark} 
In the scalar case, equation (\ref{smlapl}) is just the $p$-Laplacian. The
case of differential forms on the manifold $M=\r^n$ appears  in
section 6.1 of \cite{OI} where it is investigated by the method of Hodge dual systems, see also \cite[\S 8]{IL93}.

\medskip

\begin{thm} \label{th.solveNLeq}%
Assume $T_{q,p}^{k}(M)=0$, $(1<q,p<\infty)$  and $\alpha\in
L^{q'}(M,\Lambda^{k})$ where $q'=q/(q-1)$.

(A) If $\ds \int_M\langle\alpha,\varphi\rangle \, d\mathrm{vol}=0$
 for any $\varphi\in Z_{q}^{k}(M)$, then  (\ref{wlapl}) has
a solution $\theta\in\Omega_{q,p}^{k}(M)$.

(B) Conversely,  if (\ref{wlapl}) is solvable in
$\Omega_{q,p}^{k}(M)$, then $\ds
\int_M\langle\alpha,\varphi\rangle \, d\mathrm{vol}=0$  for any
$\varphi\in C_c^{\infty}(M,\Lambda^{k})$  such that $d\varphi=0$.
\end{thm}

\medskip

\textbf{Proof}
Assertion (B) follows from the previous Lemma, because for any
$\varphi\in C_c^{\infty}(M,\Lambda^{k})\cap \ker d$, we have
\[
\int_M\langle\alpha,\varphi\rangle \, d\mathrm{vol} =
\int_M\langle \left\Vert d\theta\right\Vert
^{p-2}d\theta,d\varphi\rangle \, d\mathrm{vol} =0.
\]
Let us prove assertion (A). The
variational functional corresponding to (\ref{wlapl}) reads
\[
I(\theta)=\frac{1}{p}\int_{M}\left\Vert d\theta\right\Vert ^{p}
\, d\mathrm{vol}-\int_M\langle\alpha,\theta\rangle  \,
d\mathrm{vol}.
\]
We first show that the functional
$I(\theta):\Omega^k_{q,p}(M)\to\r$ is bounded from below:

For any $\theta\in \Omega^k_{q,p}(M)$ there exists a unique element
$z_{q}(\theta)\in Z_{q}^{k}(M)$ such that $\left\Vert
\theta-z_{q}(\theta)\right\Vert _{q}\leq\inf_{z\in
Z_{q}^{k}(M)}\left\Vert \theta-z\right\Vert _{q}$; this follows
from the uniform convexity of $\Omega^k_{q,p}(M)$. Since
$T_{q,p}^{k}(M)=0$, the Proposition \ref{th.sobin2} implies that
\begin{equation}\label{ineq.sob}
 \left\Vert \theta-z_{q}(\theta)\right\Vert _{q}\leq C\left\Vert
 d\theta\right\Vert _{p}
\end{equation}
for some positive constant $C$. Using this inequality and
H\"older's inequality, we obtain
\[
 I(\theta)\geq\frac{1}{p}\left\Vert d\theta\right\Vert
_{p}^{p}-\left\Vert \alpha\right\Vert _{q'}\left\Vert
\theta-z_q(\theta)\right\Vert _{q}\geq\frac{1}{p}\left\Vert
d\theta\right\Vert _{p}^{p}-C\left\Vert \alpha\right\Vert
_{q'}\left\Vert d\theta\right\Vert _{p}.
\]
Since the function $f : \r \to \r$ defined by
$f(x)=\frac{1}{p}|x|^{p}- ax$ is bounded below for $x\geq0$, the
previous inequality implies that
$$
 \inf_{\theta\in\Omega_{q,p}^{k}(M)} I(\theta) > -\infty.
$$

We now prove the existence of a minimizer of $I$ on
$\Omega_{q,p}^{k}(M)$: Let $\left\{
\theta_{i}\right\}\subset\Omega_{q,p}^{k}(M) $ be a sequence such
that $I(\theta_{i})\rightarrow\inf\, I(\theta)$. Because the
function $f(x)=\frac{1}{p}|x|^{p}- ax$ is proper, the inequality
\[
I(\theta_{i})\geq \frac{1}{p}\left\Vert d\theta_{i}\right\Vert
_{p}^{p}-C\left\Vert \alpha\right\Vert _{q'}\left\Vert
d\theta_{i}\right\Vert _{p}
\]
implies that $\{ \left\Vert d\theta_{i}\right\Vert _{p}\}\subset
\r$ is bounded and, by (\ref{ineq.sob}), $\{\left\Vert
\theta_{i} -z_{q}(\theta_{i})\right\Vert _{q}\}$ is also bounded. Hence
the sequence $\{\widetilde{\theta}_{i}:=
\theta_{i}-z_{q}(\theta_{i})\}$ is bounded in
$\Omega_{q,p}^{k}(M)$.

\medskip

Since $\Omega_{q,p}^{k}(M)$ is reflexive there exists a
subsequence (still noted $\{ \widetilde{\theta}_{i}\}$) which
converges weakly to some $\theta_{0}\in \Omega_{q,p}^{k}(M)$. By
the weak continuity of the functional
$\int_M\langle\alpha,\theta\rangle  \, d\mathrm{vol}$ in
$\Omega_{q,p}^{k}(M)$ we have
\begin{equation}\label{weak}
\lim_{i\rightarrow\infty}\int_{M}\langle\alpha,\widetilde{\theta}_{i}\rangle
\, d\mathrm{vol} = \int_{M}\langle\alpha,\theta_{0}\rangle \,
d\mathrm{vol}
\end{equation}
The lower semicontinuity of the norm under the weak convergence
implies that
\[
\left\Vert d\theta_{0}\right\Vert _{p}\leq
\liminf_{i\rightarrow\infty}\left\Vert
d\widetilde{\theta}_{i}\right\Vert _{p}.
\]
Combining the last inequality with (\ref{weak}) we obtain
\[
I(\theta_{0})\leq\liminf_{i\rightarrow\infty}\, I(\theta_{i})
\]
and by the choice of $\theta_{i}$ we  finally have $I(\theta_{0})
=  \inf\, I(\theta)$.

It is now clear that $\theta_{0}$ is a solution of (\ref{wlapl}),
hence a weak solution of (\ref{smlapl}).

\qed

\medskip

\textbf{Definition.} The Riemannian manifold $(M,g)$ is \emph{$s$-parabolic} if for any $\varepsilon>0$, 
there exists  a smooth function  $f_{\varepsilon}$ with compact support, such that $ f_{\varepsilon}=1$ 
on the ball $B(x_0,1/{\varepsilon})$ and $\| d f_{\varepsilon}\|_{L^s(M)} \leq \varepsilon$.
where $x_0\in M$ is a fixed base point.

\medskip

Some basic facts about this notion can be found in \cite{tr99}. 

\medskip

\begin{cor} Assume as above that $T_{q,p}^{k}(M)=0$  and
$\alpha\in L^{q'}(M,\Lambda^{k})$ where $q'=q/(q-1)$ , $(1<q,p<\infty)$.

Assume furthermore that $M$ is
$s$-parabolic for $\frac{1}{s} = \frac{1}{p}+\frac{1}{q}$.

Then equation  (\ref{wlapl}) is solvable in
$\Omega_{q,p}^{k}(M)$,  if and only if
$\int_M\langle\alpha,\varphi\rangle \, d\mathrm{vol}=0$  for any
$\varphi\in Z_{q}^{k}(M)$.
\end{cor}

\textbf{Proof} The condition is sufficient by the previous theorem.
Now let $\varphi\in Z_{q}^{k}(M)$ be arbitrary and let
$R^M_{\varepsilon}$ be the smoothing operator and
$f_{\varepsilon}$ be as in the previous definition. Then
$$\varphi_{\varepsilon} := f_{\varepsilon} R^M_{\varepsilon}(\varphi)
\in C_c^{\infty}(M,\Lambda^{k}).$$

 Let us observe that
$$
 \|  |d\theta|^{p-2}d\theta\|_{L^{p'}(M)} = \| d\theta\|^{p/p'}_{L^p(M)}
 $$
 where $p'=p/(p-1)$. Since   $\frac{1}{s} = 1-\frac{1}{p'}+\frac{1}{q}$,
we have by H\"older's inequality:
\begin{eqnarray*}
 \int_M\langle\alpha,\varphi_{\varepsilon}\rangle \, d\mathrm{vol}
&=&
\int_M\langle \left\Vert d\theta\right\Vert
^{p-2}d\theta,d\varphi_{\varepsilon}\rangle \, d\mathrm{vol}
\\ &=&
\int_M\langle \left\Vert d\theta\right\Vert
^{p-2}d\theta,d f_{\varepsilon} \wedge R^M_{\varepsilon}(\varphi)\rangle \, d\mathrm{vol}
\\ &\leq &
 \|  |d\theta|^{p-2}d\theta\|_{L^{p'}(M)}\| d f_{\varepsilon}\|_{L^s(M)}
\|R^M_{\varepsilon}(\varphi)\|_{L^q(M)}
\\ &\leq &
 \left(\| d\theta\|^{p'/p}_{L^p(M)}
\|R^M_{\varepsilon}(\varphi)\|_{L^q(M)}\right) \| d f_{\varepsilon}\|_{L^s(M)}
\end{eqnarray*}
As $\varepsilon \to 0$, we have $\| d f_{\varepsilon}\|_{L^s(M)} \to 0$ while
$ \left(\| d\theta\|^{p'/p}_{L^p(M)}
\|R^M_{\varepsilon}(\varphi)\|_{L^q(M)}\right) $ remains bounded. On the other hand,
$$\lim_{\varepsilon \rightarrow 0}
 \int_M\langle\alpha,\varphi_{\varepsilon} \rangle \, d\mathrm{vol}
 =
  \int_M\langle\alpha,\varphi \rangle \, d\mathrm{vol}
$$
and the result follows.
\qed

\section{Torsion in $L_{2}$-cohomology and the Hodge-Kodaira decomposition}

In this section, we study some connection between the torsion in $L_{2}$-cohomology and 
the Laplacian $\Delta$ acting on differential forms  on the complete Riemannian manifold  $(M,g)$.

\medskip

Recall that  $\Delta=d\delta+\delta d$ where  $\delta$ is the
formal adjoint operator to the exterior differential $d$. We look
at $\Delta$ as an unbounded operator acting on the Hilbert space
$L^{2}(M,\Lambda^k)$. In particular, all function spaces
appearing in this section are subspaces of  $L^{2}(M,\Lambda^k)$.
We denote by $\mathcal{H}_{2}^{k}(M) =L^{2}(M,\Lambda^k) \cap
\ker \Delta$  the space of  $L^2$ harmonic forms.

 \medskip

We begin with  the following result, which can be proved by
standard arguments from functional analysis:
\smallskip

\begin{thm} \label{th: Imdelta}%
For any complete Riemannian manifold  $(M,g)$, the following
conditions are equivalent:
\begin{enumerate}[(a)]
\item  $\Imm\Delta$ is a closed subspace in $ L^{2}(M,\Lambda^k)$;
\item $\Imm \Delta = \left(\mathcal{H}_{2}^{k}(M)\right)^{\bot}$;
\item There exists a bounded linear operator \
$G :  L^{2}(M,\Lambda^k) \to L^{2}(M,\Lambda^k)$ \
such that  for any $\alpha \in   L^{2}(M,\Lambda^k)$ we have
$$
\Delta \circ G \, \alpha = G\circ  \Delta  \, \alpha =\alpha - H \alpha
$$
where $H :   L^{2}(M,\Lambda^k)   \to \mathcal{H}_{2}^{k}(M)$
is the orthogonal projection onto the space of  $L^2$ harmonic forms.
\end{enumerate}
\end{thm}

\medskip

\textbf{Remark:} $G$ is called the \emph{Green operator}. It is
not difficult to check that $d\circ G = G\circ d$ and
$\delta\circ G = G\circ \delta$.

\medskip

For the convenience of the reader, we briefly explain the proof of this Theorem:

\medskip

\textbf{Proof}
(a) $\Leftrightarrow$ (b): \
Because $\Delta$ is self-adjoint, we know by standard functional analysis
(see e.g.  \cite{Brezis},  page 28) that
$\overline{\Imm \Delta} =  \left(\mathcal{H}_{2}^{k}(M)\right)^{\bot}$,

\medskip

(b) $\Rightarrow$ (c): \  This follows from the Banach Open
Mapping Theorem. More precisely, let us denote by 
\begin{equation*}
 E := \{
\omega\in L^{2}(M,\Lambda^k) \; \big| \; \omega \bot
\mathcal{H}_{2}^{k}(M) \ \mathrm{and} \  \Delta \omega \in
L^2(M,\Lambda^k)\}
\end{equation*}
 the domain of the Laplacian. This is a Hilbert space for
the graph norm $\|\omega\|_{E} := \|\omega\|_{L^2} +
\|\Delta\omega\|_{L^2}$ and the map $\Delta : E \to \Imm\Delta =
\left(\mathcal{H}_{2}^{k}(M)\right)^{\bot}$ is a continuous
bijective operator.

From the Banach Open Mapping Theorem, we know that the map
$$G := \Delta^{-1} \circ (1-H) :  L^{2}(M,\Lambda^k) \to L^{2}(M,\Lambda^k)$$
given by the composition
$$
  L^{2}(M,\Lambda^k)   \stackrel{1-H}{\longrightarrow} \left(\mathcal{H}_{2}^{k}(M)\right)^{\bot}
  \stackrel{\Delta^{-1} }{\longrightarrow}
  E \subset  L^{2}(M,\Lambda^k)
$$
is continuous. It is clear that $G$ satisfies the required properties.

\medskip

(c) $\Rightarrow$ (b): \  Condition (c)  obviously implies that ${\Imm \Delta} \supset\left(\mathcal{H}_{2}^{k}(M)\right)^{\bot}$. 
The other inclusion ${\Imm \Delta} \subset\left(\mathcal{H}_{2}^{k}(M)\right)^{\bot}$ always holds since $\Delta$ is self-adjoint.
\smallskip

\qed

\medskip

In the case of complete Riemannian manifolds, we have the
following~:

\begin{thm}\label{thm:Hodge-Kodaira} %
For any complete Riemannian manifold  $(M,g)$, we have
$$\mathcal{H}_{2}^{k}(M) =\ker d\cap\ker\delta\cap L^{2}(M,\Lambda^k),$$
and the orthogonal decomposition
$$
 L^{2}(M,\Lambda^k) =
 \overline{\Imm   d} \oplus  \overline{\Imm  \delta }\oplus\mathcal{H}_{2}^{k}(M).
 $$
\end{thm}
The first part is due to Andreotti and Vesentini, the second part
is the well known Hodge-Kodaira decomposition. A proof is given
in \cite[Theorem 24 and 26]{derham}.

\qed

\medskip

Using both previous Theorems, we can now prove the following result:

\begin{thm} \label{thm:L2torsion}%
 For any complete Riemannian manifold  $(M,g)$, the following conditions are equivalent:
\begin{enumerate}[(i.)]
\item $\Imm \Delta = \left(\mathcal{H}_{2}^{k}(M)\right)^{\bot}$;
\item we have the orthogonal decomposition
$$ L^{2}(M,\Lambda^k) =
 {\Imm  d} \oplus   {\Imm \delta }\oplus\mathcal{H}_{2}^{k}(M);$$
 \item   $\Imm d$ and  $\Imm \delta$ are closed  in $L^{2}(M,\Lambda^k)$;
 \item $T^k_2 (M) = 0$ and $T^{n-k}_2 (M) = 0$.
 \end{enumerate}
\end{thm}

\medskip

We will also need the following

\medskip

\begin{lem} If  $ T_{2}^{k}(M)=0$, then
$$ \Imm (\delta d) =  \Imm (\delta ) $$
as subsets of $L^{2}(M,\Lambda^k)$.
\end{lem}

\medskip

\textbf{Proof} It is clear that $\Imm (\delta d) \subset  \Imm
(\delta )$.  To prove the other inclusion, consider an arbitrary
element $\alpha \in \Imm \delta$. Because $\Imm \delta \bot \ker
d = Z_2^k(M)$, we know by Theorem \ref{th.solveNLeq} that we can
find a form  $\theta \in L^{2}(M,\Lambda^k)$ such that $\delta d
\, \theta = \alpha$. In particular  $\alpha \in \Imm \delta d$.

\qed

\medskip

\textbf{Remark. }  Using  the formula $\delta = \pm *d*$, we see
that this lemma also says that $\Imm (d\delta ) =  \Imm (d)$,
provided $T_{2,2}^{n-k}(M)=0$.

 \newpage

\textbf{Proof of Theorem \ref{thm:L2torsion}.} 

(i) $\Rightarrow$ (ii): \  Condition (i) is equivalent to (c) of
Theorem \ref{th: Imdelta}. Hence, assuming (i), we know that  any
$\alpha \in  L^{2}(M,\Lambda^k)$ can be written as
$$\alpha - H \alpha = \Delta \circ G \, \alpha
= d(\delta G \alpha) + \delta (d G \alpha) $$
and the decomposition (ii) follows.

\smallskip

(ii) $\Rightarrow$ (iii): \  is clear from Theorem \ref{thm:Hodge-Kodaira}.

\smallskip

(iii) $\Leftrightarrow$ (vi): \ Follows from the definition of
torsion and the formula $\delta = \pm *d*$.

(vi) $\Rightarrow$ (i): \  We know from the  previous lemma and the orthogonality of $\Imm d$ and $\Imm \delta$ that
$$\Imm \Delta = \Imm (d\delta + \delta d) = \Imm (d\delta) + \Imm (\delta d)
=  \Imm (d) + \Imm (\delta),$$
provided $T^k_2 (M) = T^{n-k}_2 (M) = 0$. In particular, $\Imm \Delta$ is closed, since
$\Imm d$ and  $\Imm \delta$ are closed, and we conclude by Theorem \ref {th: Imdelta}
that  $\Imm \Delta = \left(\mathcal{H}_{2}^{k}(M)\right)^{\bot}$.

\qed

\medskip

\begin{cor} If $(M,g)$ is complete, then the equation $\Delta\omega=\alpha \in  L^{2}(M,\Lambda^k)$ is solvable in
 $L^{2}(M,\Lambda^k)$ for any  $\alpha\bot\mathcal{H}_{2}^{k}(M)$, if and only if
 $$
  T_{2}^{k}(M)=0 \quad  \mathrm{and } \quad T_{2}^{n-k}(M)=0.
 $$
\end{cor}

The proof is immediate.

\qed

\medskip

In conclusion, we formulate the following version of Hodge
Theorem and Poincar\'{e} duality for $L^2$-cohomology:

\begin{cor} If $(M,g)$ is a complete Riemannian manifold such that
\\
$T_{2}^{k}(M)=T_{2}^{n-k}(M)=0$, then
$$
\overline{H}^k_2(M) =  H^k_2(M) \cong \mathcal{H}_{2}^{k}(M) \cong
 \mathcal{H}_{2}^{n-k}(M) \cong H^{n-k}_2(M) =   \overline{H}^{n-k}_2(M).
$$
\end{cor}

\medskip

\textbf{Proof} The equality $\overline{H}^k_2(M) =  H^k_2(M)$ is equivalent to $T_{2}^{k}(M)=0$.

From Theorem \ref{thm:L2torsion}, we know that if the torsion vanishes, then
$$\ker d = (\Imm \delta)^{\bot}= \Imm d \oplus  \mathcal{H}_{2}^{k}(M),$$
i.e. $H^k_2(M) \cong \mathcal{H}_{2}^{k}(M)$ by definition of cohomology.

\smallskip

The isomorphism $\mathcal{H}_{2}^{k}(M) \cong
\mathcal{H}_{2}^{n-k}(M)$ is given by the Hodge $*$ operator and
the proof now ends as it begins.

\qed

\section*{Appendix: A ``classic'' proof of Theorem \ref{th.sobin1}   in the compact case.}

In this appendix, we  shortly give another proof of Theorem
\ref{th.sobin1} for compact manifolds which is based on the Hodge
De-Rham theory and  the regularity theory for elliptic systems,
together with some techniques from functional analysis. All these
tools were available 40 years ago, however, we did not find a written proof in the literature.

\medskip

We start with the fact that the space of harmonic currents
on a compact Riemannian manifold $(M,g)$ is finite dimensional and
that we can construct two linear operators acting on currents on $M$
 $$ G,H : \mathcal{D}'(M) \to  \mathcal{D}'(M),$$
and such that
\begin{enumerate}[i)]
\item $ \ker \Delta = \Imm H = \ker (I-H)$;
\item $ \ker \Delta \cap \Imm(I-H) = \{0\}$;
\item $\Delta \circ G = (I-H)$;
\item $\Delta \circ (I-H) = \Delta$;
\item $d\circ G = G\circ d$.
\end{enumerate}
This result is theorem 23 in \cite{derham}, the operator $H$ is
the projection onto the space of harmonic forms and $G$ is the
Green operator.

\medskip

Using elliptic regularity, we can prove the following theorem:

\begin{thm}\label{th.regG} %
The Green operator defines a bounded linear operator
 $$G : W^{m,p}(M,\Lambda^k) \to W^{m+2,p}(M,\Lambda^k)$$
for any $m\in \mathbb{N}$. Here $W^{m,p}(M,\Lambda^k)$ is the
Sobolev space of differential forms of degree $k$ on $M$ with
coefficients in $W^{m,p}$.
\end{thm}

\medskip

Assuming this result for the time being, let us conclude the proof of
Theorem \ref{th.sobin1}. We first state the following corollary:

\begin{cor}
For any compact Riemannian manifold $(M,g)$, there exists a constant $C_1$ such that
\begin{equation}\label{ineq.C1}
   \|\theta - \zeta \|_{W^{1,p}(M)} \leq C_1 \|d\theta \|_{L^{p}(M)},
\end{equation}
 where \
 $\zeta := H \, \theta + d\delta G \, \theta$.
\end{cor}

\medskip

\textbf{Proof}  From previous theorem, we see that
$\delta \circ G : L^p (M,\Lambda^k) \to W^{1,p}(M,\Lambda^{k+1})$ is a bounded operator.

\medskip

Since $\Delta G = (d\delta + \delta d) G = (I-H)$, we have
$\theta - \zeta = \delta d G \, \theta = \delta  G \, d\theta$ and thus
 $$\|\theta - \zeta \|_{W^{1,p}(M)} =
 \| \delta  G \, d\theta \|_{W^{1,p}(M)}
 \leq C_1 \|d\theta \|_{L^{p}(M)},$$
where $C_1$ is the operator norm
$C_1     :=  \| \delta  G  \|_{L^p \to W^{1,p}}.$

\qed

\medskip

\subsection*{Proof of Theorem \ref{th.sobin1}.}
The classical Sobolev embedding theorem on compact manifolds,
states in particular that there is a constant $C_2$ such that
\begin{equation}\label{ineq.C2}
  \|\omega\|_{L^{q}(M)} \leq C_2 \|\omega\|_{W^{1,p}(M)},
\end{equation}
provided that conditions (\ref{cond.sobolev}), are satisfied.

\smallskip

Combining (\ref{ineq.C1}) and (\ref{ineq.C2}) and observing that,
by the Sobolev embedding theorem and (\ref{cond.sobolev}), we have
$\zeta = H \, \theta + d\delta G \, \theta \in Z^k_q(M)$, we
obtain (\ref{ineq.sob1}) with $C=C_1C_2$.

\qed

\bigskip

\subsection*{Proof of Theorem \ref{th.regG}}
The proof is in several steps.

\smallskip

Step 1. The elliptic estimate for the Laplacian acting on forms
on a compact manifold says that there exists a constant $A_m$
such that for any form $\theta \in W^{m+2,p}(M,\Lambda^k)$ we have
\begin{equation}\label{ineq.ellipticestimate}
   \|\theta \|_{W^{m+2,p}(M)} \leq A_m \left(\|\Delta \theta \|_{W^{m,p}(M)}
   +  \|\theta \|_{W^{m,p}(M)}\right).
\end{equation}
This result is deep. The case $p=2$ is proved in proved in  \cite[\S 6.29]{warner}, 
the scalar case for any $p\in (0,\infty)$ can be found in  \cite[\S 9.5]{GT} and the general 
case in  \cite[Chapter IV]{ADN2}.

\medskip

Step 2. A first consequence of this estimates is the
hypoellipticity of the Laplacian, i.e. the fact if $\Delta
\theta$ is a smooth form, then $\theta$ itself is smooth (the
proof follows from a bootstrap argument based on
(\ref{ineq.ellipticestimate}) and the fact that $\cap_{m\geq
1}W^{m,p}(M) = C^{\infty}(M)$.) It follows in particular that the
Green operator $G$ maps smooth forms to smooth forms.

\medskip

Step 3.  Using (\ref{ineq.ellipticestimate}), we show that for
any sequence $\{\theta_i\} \subset W^{m+2,p}$, we have
\begin{equation}\label{cond.bornee}
  \|\Delta \theta_i \|_{W^{m,p}(M)} \quad \mathrm{bounded} \
  \Rightarrow \ \|(I-H)\theta_i \|_{W^{m,p}(M)} \quad
  \mathrm{bounded}.
\end{equation}
Indeed, otherwise there exists a sequence such $\|\Delta \theta_i
\|_{W^{m,p}(M)}$ is bounded and $\|(I-H)\theta_i \|_{W^{m,p}(M)}
\to \infty$. Let us set
$$
\varphi_i := \frac{(I-H)\theta_i }{\|(I-H)\theta_i
\|_{W^{m,p}(M)}} \in W^{m+2,p}(M),
$$
we then have $ \|\varphi_i \|_{W^{m,p}(M)} = 1$ and
$$
\lim_{i\to \infty} \|\Delta \varphi_i \|_{W^{m,p}(M)} =
\frac{\|\Delta \theta_i \|_{W^{m,p}(M)}}{\|(I-H)\theta_i
\|_{W^{m,p}(M)}} = 0.
$$

The elliptic estimate (\ref{ineq.ellipticestimate}) gives us
$$
   \|\varphi_i \|_{W^{m+2,p}(M)} \leq A_m \left(\|\Delta \varphi_i\|_{W^{m,p}(M)}
   +  \|\varphi_i\|_{W^{m,p}(M)}\right)
$$
and thus $\{\varphi_i \}$ is bounded in $W^{m+2,p}(M)$.

\smallskip

Because $W^{m+2,p}(M)$ is reflexive, there exists a subsequence
which converges weakly in $W^{m+2,p}(M)$. We still denote this
subsequence  by $\{\varphi_i \}$. Let $\varphi \in W^{m+2,p}(M)$
be the weak limit of this subsequence, we then have by the lower
semi-continuity of the norm
$$
 \|\Delta \varphi\|_{W^{m,p}(M)} \leq \liminf_{i\to \infty} \|\Delta
 \varphi_i\|_{W^{m,p}(M)}
  = 0,
$$
hence $\varphi \in \ker \Delta$. Since we also have $\varphi \in
\Imm(I-H)$ we must have $\varphi =0$.

\smallskip

By the compactness of the embedding  $W^{m+2,p}(M) \subset
W^{m,p}(M)$, we may assume that this subsequence converges
strongly in $W^{m,p}(M)$. In particular we have
$$
 1 = \lim_{i \to \infty}\|\varphi_i\|_{W^{m,p}(M)}
 =  \|\lim_{i \to \infty}\varphi_i\|_{W^{m,p}(M)}
 = 0,
$$
This contradiction proves (\ref{cond.bornee}).

\medskip

Step 4. We now show that:
$$
  \Delta \left( W^{m+2,p}(M)\right) \quad \text{is closed in} \quad
  W^{m,p}(M)
$$
Indeed, for any $\omega \in W^{m,p}(M)$ in the closure of $\Delta
\left( W^{m+2,p}\right)$, there exists a sequence $\{\theta_i \}
\subset W^{m+2,p}$, such that  $\Delta\theta_i \to \omega$. By
step 3,  $\{(I-H)\theta_i\}$ is bounded in $W^{m,p}$, and by
(\ref{ineq.ellipticestimate}), this sequence is also bounded in
$W^{m+2,p}$ (recall that $\Delta (I-H)\theta_i= \Delta\theta_i$).

\smallskip

By the compactness of the embedding  $W^{m+2,p}(M) \subset
W^{m,p}(M)$, there exists a subsequence such that
$\{(I-H)\theta_i\}$ converges strongly in $W^{m,p}$, and by
(\ref{ineq.ellipticestimate}) again, $\{(I-H)\theta_i\}$ converges
in $W^{m+2,p}$.

Let us denote by $\ds \psi = \lim_{i\to \infty} (1-H)\theta_i$, we
then have $\omega = \Delta \psi \in \Delta \left(
W^{m+2,p}(M)\right)$.

\medskip

Step 5.  Let us denote by $\mathcal{E}^{m,p}= \ker H \cap
W^{m,p}(M,\Lambda^k) = \Imm (I-H) \cap W^{m,p}(M,\Lambda^k) $.
Then $\Delta : \mathcal{E}^{m+2,p}\to \mathcal{E}^{m,p}$ is
continuous, injective and has closed image by previous step.
Furthermore, $\Imm \Delta \subset  \mathcal{E}^{m,p}$ is dense because any
smooth form in $ \mathcal{E}^{m,p}$ is the image under $\Delta$
of a smooth form in $\mathcal{E}^{m+2,p}$. To sum up, we have
proved that
$$\Delta : \mathcal{E}^{m+2,p}\to \mathcal{E}^{m,p}$$
is a continuous linear bijection.

\medskip

Step 6.  By the Banach open mapping theorem, we finally  see that
 $$G = \Delta^{-1} \circ (1-H) :  W^{m,p}(M,\Lambda^k) \to
 \mathcal{E}^{m+2,p} \subset W^{m+2,p}(M,\Lambda^k)$$
 is a bounded operator.

 \qed

\end{document}